\documentclass[12pt]{article}    % Specifies the document style.
\usepackage{setspace}
\usepackage{hyperref}
\usepackage{graphicx}
\usepackage{epsfig}
\usepackage{url}
\usepackage{amsmath}
\usepackage{mathrsfs}
\usepackage{setspace}
\usepackage{hyperref}

\usepackage{epsfig}
\usepackage{url}
\usepackage{amsmath}
\usepackage{mathrsfs}
\usepackage{amssymb}

\usepackage{amsmath}
\usepackage{bm}
\usepackage{amsfonts}
\usepackage{amsmath}
\usepackage{mathrsfs}
\usepackage{multirow}

\usepackage{tablefootnote}
\usepackage{booktabs}
\usepackage{multirow}
\usepackage{siunitx}
\usepackage[flushleft]{threeparttable}
\usepackage{color}
\usepackage{hyperref}
\usepackage[caption = false]{subfig}

\usepackage{titlesec}
\usepackage{listings}

\usepackage{amssymb}

\usepackage{amssymb}

\numberwithin{equation}{section}

\def \R {{\rm I\kern -2.2pt R\hskip 1pt}}

\newcommand{\qed}{\rule{0.5em}{1.5ex}}

\newtheorem{remark}{Remark}

\begin{document}
	
	\begin{center} { \large \sc  Bivariate distributions with equi-dispersed normal conditionals and related models}\\
		\vskip 0.1in {\bf Barry C. Arnold}\footnote{barry.arnold@ucr.edu}
		\\
		\vskip 0.1in 
		Department of Statistics, University of California, Riverside, USA. \\
		\vskip 0.1in {\bf  B.G. Manjunath }\footnote{bgmanjunath@gmail.com}
		\\
		\vskip 0.1in 
		School of Mathematics and Statistics, University of Hyderabad, Hyderabad, India. \\

			\bigskip
		
		4 September, 2022

\end{center}
	 
	\begin{abstract}
	A random variable is equi-dispersed if its mean equals its variance. A Poisson distribution is a classical example of this phenomenon. However, a less well-known fact is that the class of normal densities that are equi-dispersed constitutes a one parameter exponential family.
	In the present article our main focus is on univariate and bivariate models with equi-dispersed normal component distributions. We discuss distributional features of such models, explore inferential aspects and  include an example of application of equi-dispersed models. Some related models are discused in Appendices.
	
\end{abstract}

	
	\bigskip

	\section{Introduction }  Conditionally specified bivariate models often provide useful flexible models exhibiting a variety of dependence structures. Probably the first such model to appear in the literature was the normal conditionals distribution first discussed, though not christened, by Bhattacharyya (1943). The model was reconsidered by Castillo and Galambos (1987), but perhaps the most extensive treatment of the model may be found in Arnold, Castillo and Sarabia (1999). We will summarize briefly the properties and characterization of the normal conditionals model. However our chief focus is on univariate and bivariate models with what we call equi-dispersed normal component distributions. We will say that a random variable is equi-dispersed if its mean equals its variance. For example, Poisson distributions provide well-known examples of this phenomenon. But equi-dispersion is very commom. Consider any random variable, $X$ whose
	positive mean is not equal to its variance, for definitenes suppose that $var(X)=kE(X)$. There then exists a positive multiple of $X$ that is equi-dispersed, namely $Y=X/k.$ The class of univariatre normal distributions forms a two parameter exponential family, as is well-known. Perhaps less well-known (outside of exercises in texts dealing with exponential families), is the fact that the class of normal densities which are equi-dispersed also forms an exponential family, a one parameter family in this case. 
	
	In this paper we will consider the class of bivariate distributions with equi-dispersed normal conditionals.
	Using the result  in Arnold and Strauss (1991), we know that this will constitute  a three parameter exponential family of bivariate densities. Rather than apply the Arnold-Strauss result, we will approach the problem by putting constraints on the (Bhattacharyya) class of distributions with normal conditionals. We will use the same approach to investigate the class of bivariate densities with conditional variances equal to squared conditional means, a setting in which the Arnold-Strauss approach is not possible. As more flexible alteratives to the conditionally specified models considered, we suggest that certain pseudo models (in the Filus-Filus sense, see for example Filus, Filus and Arnold (2009)) might merit consideration. We begin by reviewing the equi-dispersed normal model and its related bivariate extensions.
	
	\section{Equi-dispersed normal distributions, univariate and bivariate}
	
	We will say that a random variable $X$ has an equi-dispersed normal distribution if it has a normal distribution with its variance equal to its mean, i.e., if $X\sim Normal(\tau,\tau)$ for some $\tau \in (0,\infty).$ The density of such a random variable $X$ is of the form
	\begin{eqnarray} \label{eq-dis-norm}
	f_X(x;\tau)&=&\frac{1}{\sqrt{2\pi\tau}}e^{-(x-\tau)^2/2\tau} \nonumber \\
	\\
	&=&\frac{1}{\sqrt{2\pi\tau}}e^{-\tau/2}e^xe^{-x^2/2\tau}, \nonumber
	\end{eqnarray}
which is clearly an exponential family, and a sample of size $n$ from this distribution will have sufficient statistic $\sum_{i=1}^n X_i^2.$

Since equi-dispersion is a sub-model of the classical normal model, it is natural to test for its applicability before using the restricted model to analyze data. A standard testing procedure is available, and is described in the following sub-section.
	\subsection{Likelihood ratio test for the  univariate equi-dispersed normal distribution}
	We know that, the general form of a generalized likelihood ratio test statistic is as follows
	
	\begin{eqnarray}
	\Lambda = \frac{\sup_{\theta \in \Theta_0}L(\theta)}{\sup_{\theta \in \Theta} L(\theta)}
	\end{eqnarray}
	Here, $\Theta_0$ is a subset of $\Theta$, $L(\theta)$ is a likelihood function for the given data and we are envisioning testing $H_0: \theta \in \Theta_0$.  We reject the null hypothesis for small values of $\Lambda$.
	
	Let $X_1,...,X_n$ be a random sample from a normal distribution with mean $\mu$ and variance  $\sigma^2$. In the following we construct a likelihood ratio test  for testing $H_0: \mu=\sigma^2= \tau$.  The natural parameter space for the unrestricted model is $\Theta =\{ (\mu, \sigma^2): -\infty < \mu < \infty, \sigma^2>0 \}$.While, under the null hypothesis the parameter space is $\Theta_0 = \{ \tau=\mu=\sigma^2: \tau >0\}$. We know that maximum likelihood estimators of $\mu$ and $\sigma$ are
	\begin{eqnarray*}
	\hat{\mu} &=& \frac{1}{n} \sum_{i=1}^{n} X_i  \\
		\hat{\sigma^2} & = & \frac{1}{n} \sum_{i=1}^{n }(X_i - \hat{\mu})^2.
	\end{eqnarray*}
	
	Under $H_0$, the likelihood equation will be
	\begin{equation*}
	(d/d\tau) \ell(\tau)=-\frac{n}{2\tau}-\frac{n}{2}+ \frac{\sum_{i=1}^n X_i^2}{2\tau^2}=0,
	\end{equation*}
	which is equivalent to the equation
	\begin{eqnarray*}
		  \tau^2 + \tau - \frac{1}{n} \sum_{i=1}^{n} X^2_i =0.
	\end{eqnarray*}
    The unique positive solution to the above quadratic equation will be the m.l.e estimator of $\tau$, i.e.,
    \begin{eqnarray}
    	\hat{\tau} = \sqrt{\frac{1}{n} \sum_{i=1}^{n}X^2_i + \frac{1}{4}}+ \frac{1}{2}. 
    \end{eqnarray}

    	Therefore, the likelihood ratio test statistic will be
    \begin{eqnarray}
    \Lambda  = \Bigg(\frac{ \hat{\sigma^2}}{\hat{\tau}}\Bigg)^\frac{n}{2}  exp \biggl\{ - \frac{1}{2\hat{\tau}}  \sum_{i=1}^{n} (x_i - \hat{\tau})^2\biggr\}.
    \end{eqnarray}
    If $n$ is large, then $-2 \log \Lambda $ may be compared with a suitable $\chi^2_1$ percentile in order to decided whether $H_0$ should be accepted. 
    
    \subsection{Bivariate densities with equi-dispersed normal conditional distributions}
    	We will be interested in bivariate densities that have conditional distributions in the equi-dispersed normal family. Specifically, we consider a distribution of $(X,Y)$ with the property that, for each $y \in (-\infty,\infty)$ we have 
	\begin{equation}\label{eq-dis-con-X|Y}
	X|Y=y \sim Normal(\tau_1(y),\tau_1(y)) \mbox{   for some  }\tau_1(y) \mbox{   which may depend on  } y,
	\end{equation}
	and  for each $x \in (-\infty,\infty)$ we have 
	\begin{equation} \label{eq-dis-con-Y|X}
	Y|X=x \sim Normal(\tau_2(x),\tau_2(x)) \mbox{   for some  }\tau_2(x) \mbox{   which may depend on  } x.
	\end{equation}

	A result of Arnold and Strauss (1991), dealing with distributions with conditionals in exponetial families, may be applied here to conclude that the family of such bivariate distributions will constitute a $3$-parameter exponential family with sufficient statistics (based on a sample of size $n$) given by $$\left(\sum_{i=1}^n X_i^2,\sum_{i=1}^n Y_i^2,\sum_{i=1}^n X_i^2 Y_i^2 \right).$$
	At this point, we could refer to the Arnold and Strauss paper to identify the form of the joint density of $(X,Y)$ satisfying (\ref{eq-dis-con-X|Y}) and (\ref{eq-dis-con-Y|X}). However we will obtain this density instead by specializing in the general expression for distributions with normal conditionals introduced in Bhattacharyya (1943), using notation similar to that used in Arnold, Castillo and Sarabia (1999, p.58). If $(X,Y)$ has normal conditionals then its joint density will be of the form
	\begin{equation}\label{e3.26}
	f_{X,Y}(x,y) = \exp \left\{- \left( 1,x,x^2 \right)\left( \begin{array}{ccc} a_{00} & a_{01} &
	a_{02} \\ a_{10} & a_{11} & a_{12} \\ a_{20} & a_{21} &
	a_{22}\end{array}\right) \left( \begin{array}{ccc}1 \\ y
	\\ y^2 \end{array}
	\right)\right\}.
	\end{equation}
	with conditional moments of the form
	\begin{eqnarray}
	E(X\mid Y=y) & = & \mu _{1}(y)=-{{a_{12}y^2+a_{11}y+a_{10}} \over
		{2(a_{22}y^2+a_{21}y+a_{20})}}\label{e3.100} \\ 
	var(X\mid Y=y) & = & \sigma_{1}^2(y)={1 \over {2(a_{22}y^2+a_{21}y+a_{20})}}\label{e3.101}\\
	E(Y\mid X=x) & = & \mu _{2}(x)=-{{a_{21}x^2+a_{11}x+a_{01}} \over
		{2(a_{22}x^2+a_{12}x+a_{02})}}\label{e3.102}\\
	var(Y\mid X=x) & = & \sigma_{2}^2(x)={1 \over {2(a_{22}x^2+a_{12}x+a_{02})}}\label{e3.103}.
	\end{eqnarray}
	In order to guarantee that the marginals of (\ref{e3.26}) are non-
	negative (or equivalently to guarantee that for each fixed
	$x$, $f_{X,Y}(x,y)$ is integrable with respect to $y$ and for each fixed $y$
	it is integrable with respect to $x$), the coefficients in (\ref{e3.26}) must
	satisfy one of the two sets of conditions.
	\begin{equation}\label{e3.80}
	a_{22}=a_{12}=a_{21}=0;\;\; a_{20}>0; \;\;a_{02}>0.
	\end{equation}
	\begin{equation}\label{e3.81}
	a_{22}>0; \;\;4a_{22}a_{02}>a_{12}^2;\;\; 4a_{20}a_{22}>a_{21}^2.
	\end{equation}
	
	If (\ref{e3.80}) holds then we need to assume in addition that
	\begin{equation}\label{e3.28}
	a_{11}^2>4a_{02}a_{20}.
	\end{equation}
	in order to guarantee that (\ref{e3.26}) is
	integrable. Note that (\ref{e3.80}) and (\ref{e3.28}) yield the classical bivariate normal model.

	From these expressions for the conditional moments, it is evident that necessary and sufficient conditions for equi-dispersion of the conditional densities are that
	\begin{equation}\label{nsc-for-equi-disp}
	a_{11}=a_{12}=a_{21}=0  \ \  \ \mbox {     and      } \ \   \ a_{10}=a_{01}=-1.
	\end{equation}
	Since the normal conditionals model (\ref{e3.26}) had an 8 dimensionsal parameter space (note that $a_{00}$ is a function of the other $a_{ij}$'s chosen to normalize density to integrate to $1$) , the five constraints in (\ref{nsc-for-equi-disp}) reduce the model to a three parameter model (as expected from the Arnold and Starauss theorem). To eliminate no longer needed sub-scripts, we will relabel the three remaing parameters as
	\begin{equation}\label{new-param}
	\alpha=a_{20},  \hspace{0.5in}  \beta=a_{02} \hspace{0.5in}  \mbox{and} \hspace{0.5in} \gamma=a_{22}.
	\end{equation}	
	The equi-dispersed normal conditionals density is thus of the form
	\begin{equation}\label{eq-di-no-co}
	f_{X,Y}(x,y:\alpha,\beta,\gamma)\propto exp\{-[\alpha x^2+\beta y^2 +\gamma x^2y^2-x-y]\}
	\end{equation}
	with conditional moments

	\begin{figure}[ht!]
		\centering
		\includegraphics[width=100mm]{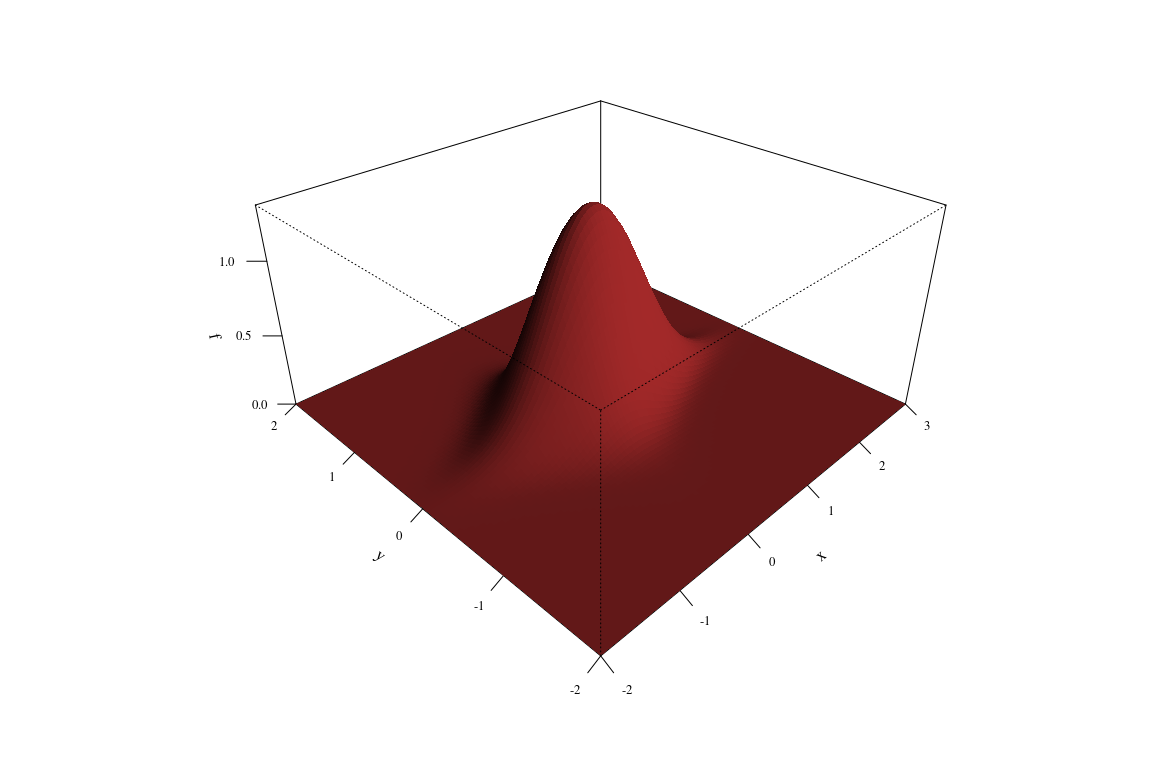}
		\caption{Density plot (strong dependence): $\alpha=1, \beta=4, \gamma=5$}
	\end{figure}
	
		\begin{figure}[ht!]
		\centering
		\includegraphics[width=100mm]{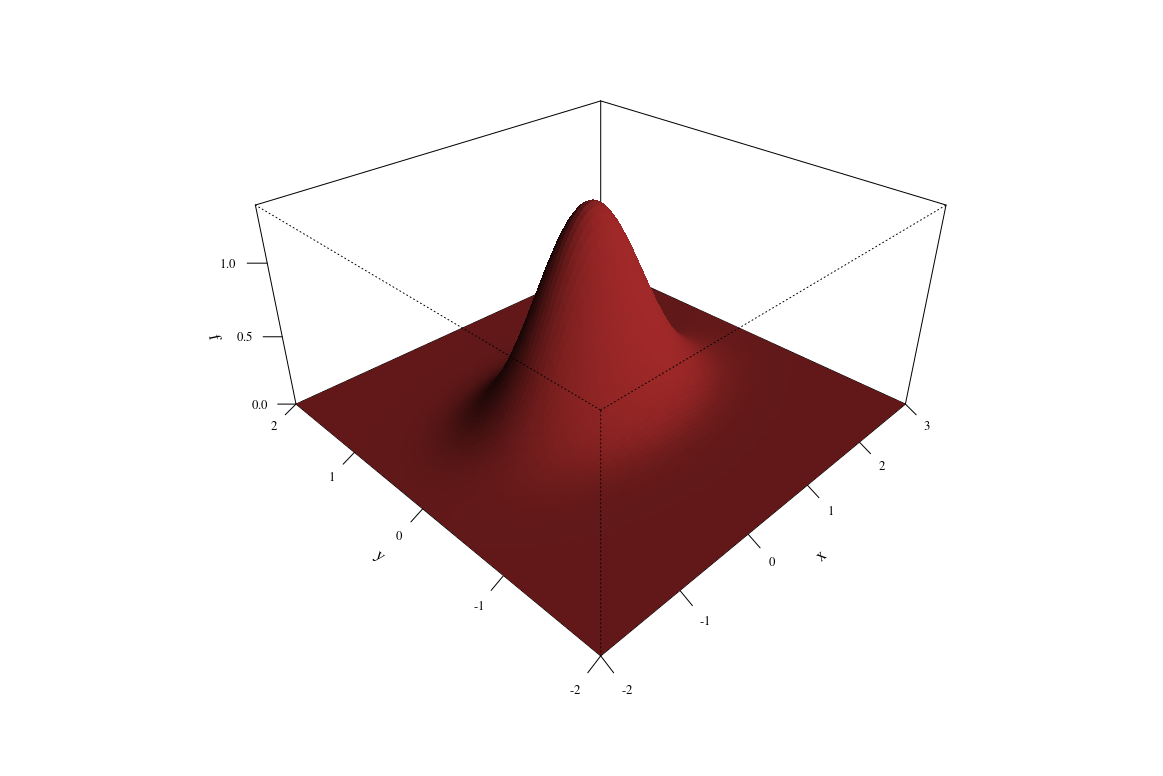}
		\caption{Density plot (near independence): $\alpha=1, \beta=4, \gamma=0.12$}
	\end{figure}
	
		\begin{figure}[ht!]
		\centering
		\includegraphics[width=100mm]{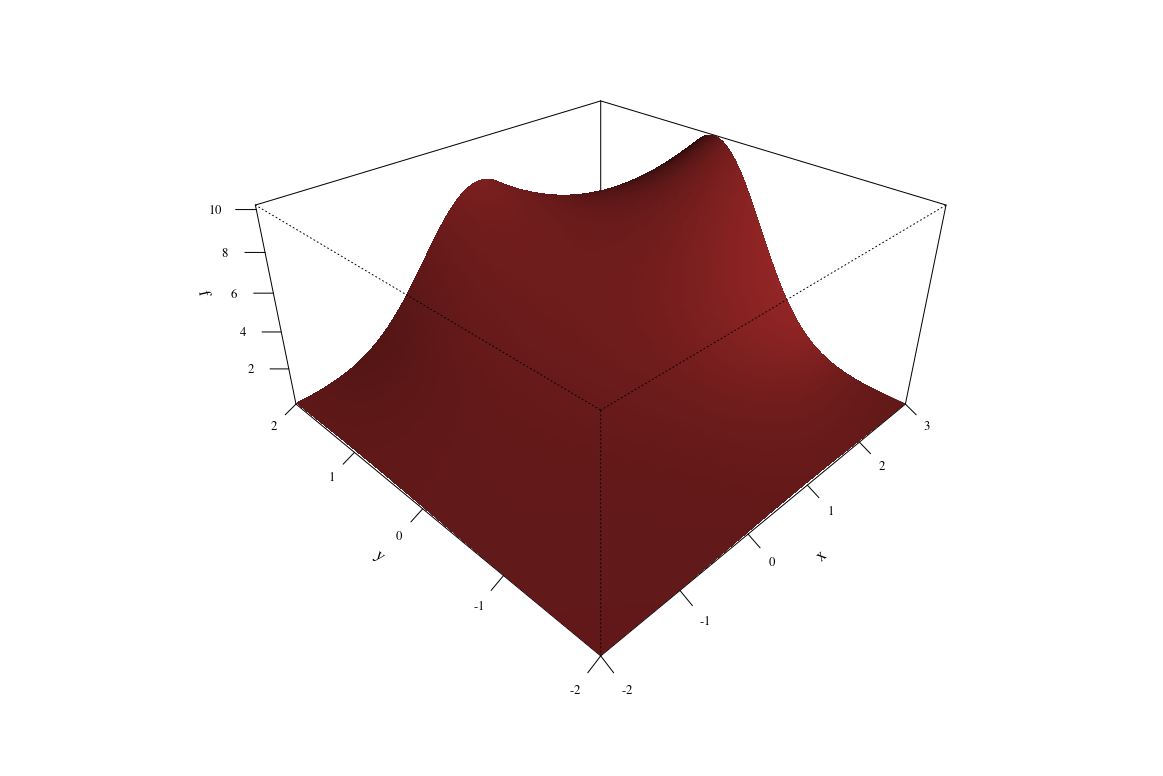}
		\caption{Density plot (bimodality): $\alpha=0.099, \beta=0.088, \gamma=0.12$}
	\end{figure}

		\begin{figure}[ht!]
		\centering
		\includegraphics[width=100mm]{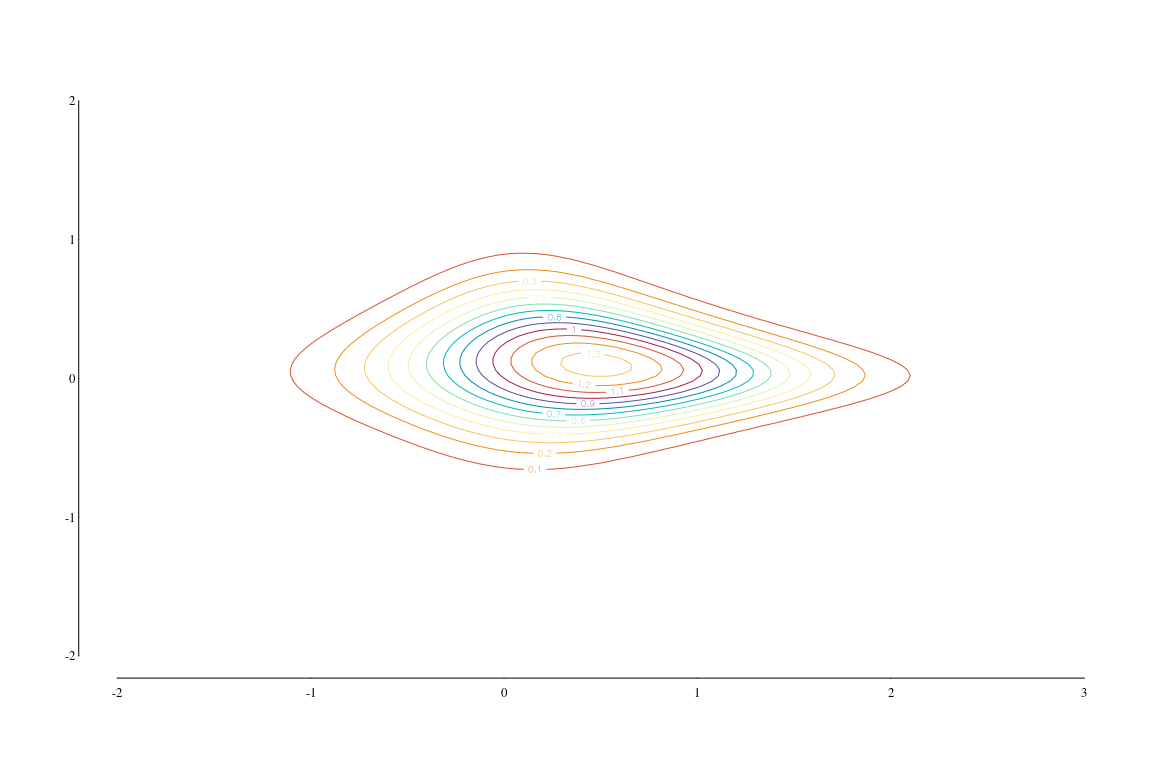}
		\caption{Contour plot (strong dependence): $\alpha=1, \beta=4, \gamma=5$}
	\end{figure}
	
	\begin{figure}[ht!]
		\centering
		\includegraphics[width=100mm]{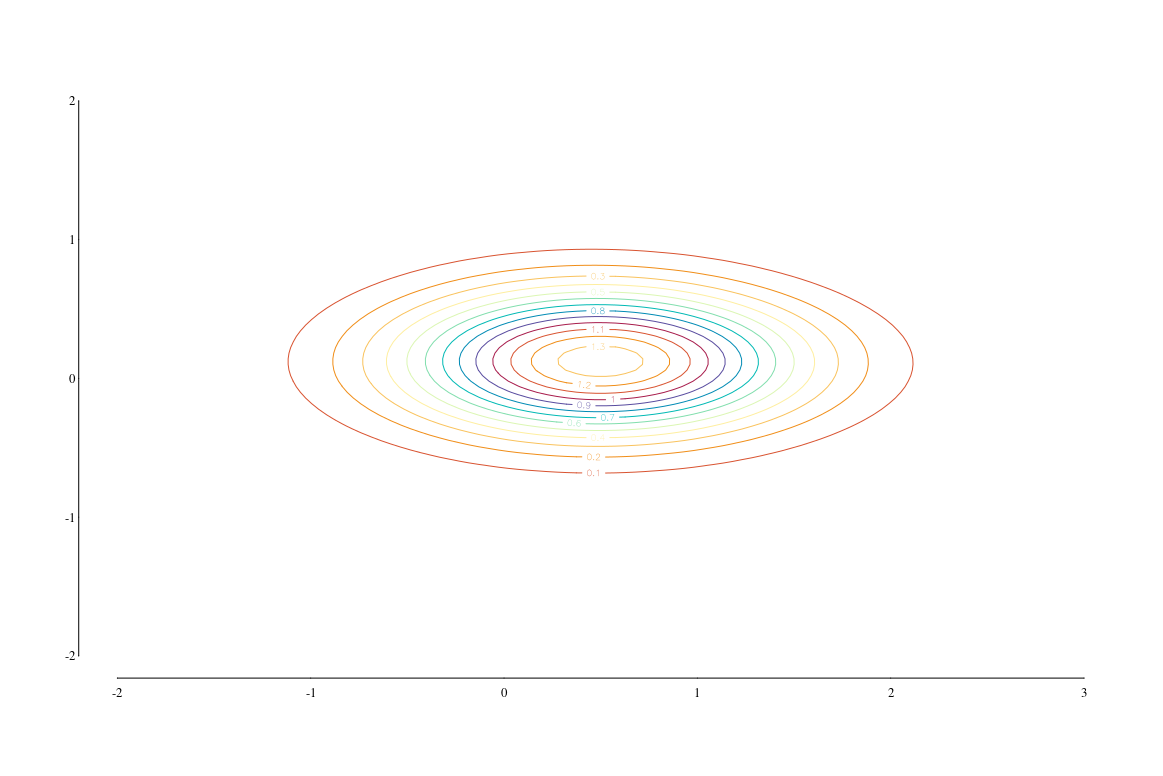}
		\caption{Contour plot (near independence): $\alpha=1, \beta=4, \gamma=0.12$}
	\end{figure}
	
	\begin{figure}[ht!]
		\centering
		\includegraphics[width=100mm]{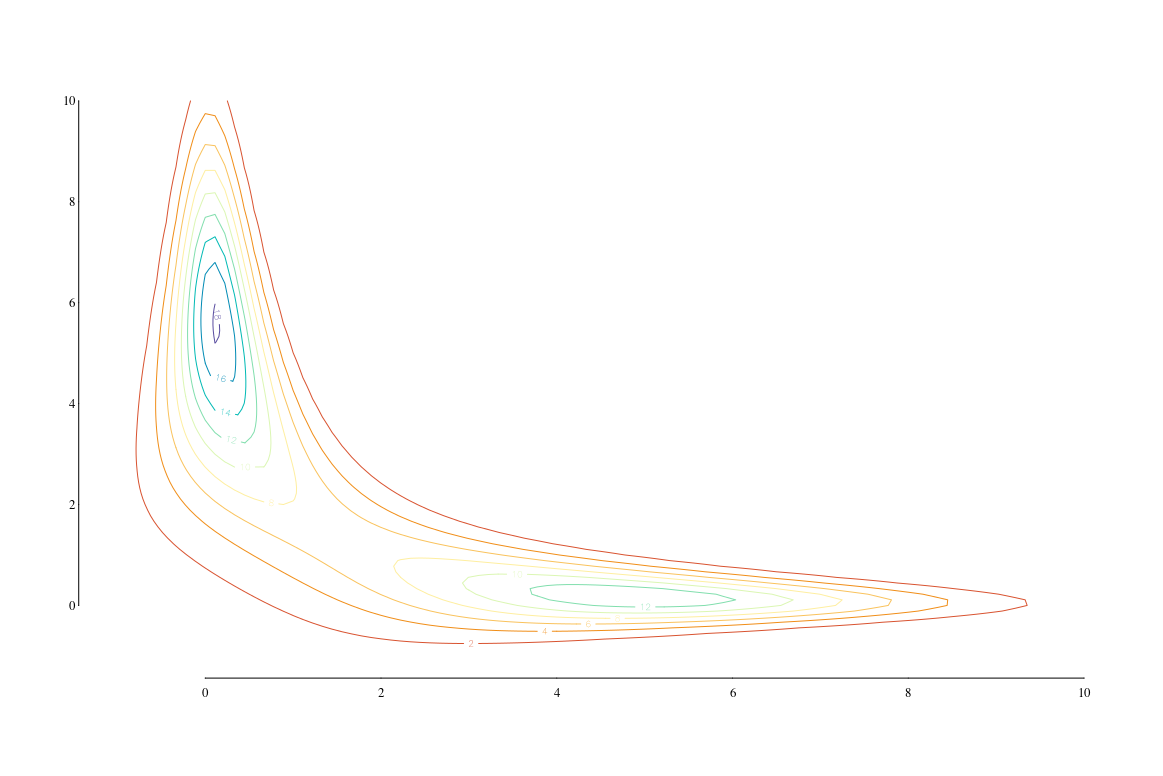}
		\caption{Contour plot (bimodality): $\alpha=0.099, \beta=0.088, \gamma=0.12$}
	\end{figure}

	\begin{eqnarray}
	E(X\mid Y=y) = var(X\mid Y=y)=\frac{1}{2(\gamma y^2+\alpha)}\label{e3.100a} \\ 
	E(Y\mid X=x) = var(Y\mid X=x)  = \frac {1}{2(\gamma x^2+\beta)}\label{e3.102a}
	\end{eqnarray}
	
	In this model we require that $\alpha>0, \beta>0$ and $\gamma \geq 0$. Note that, if $\gamma=0$, then $X$ and $Y$ are independent equi-dispersed normal variables.
	
	It is known that the full normal conditionals density (\ref{e3.26} ) can have more than one mode, (see Arnold et al. (2000) for detailed discussion of this phenomenon), although a single mode is more commonly encountered. An analogous situation is found in the case of the equi-dispersed normal conditionals density
(\ref{eq-di-no-co}). More than one mode can occur, although this is atypical. We refer to Figures  1,2,3 and Figure  4,5,6  for density and contour plots of the 
equi-dispersed normal conditionals models for different choices of parameters, exhibiting strong dependence, near independence and bimodality, respectively.

	The marginal densities are of the form
	\begin{equation}\label{X-marg}\begin{array}{l}
	f_X(x)  =  (2(\gamma x^2+\beta))^{-\frac{1}{2}}\times\\
	\hfill \hspace{0.8cm}\exp\left\{-\frac{1}{2}\left[2(\alpha x^2-x+a_{00})
	- \displaystyle\frac{1}{2(\gamma x^2+\beta)}\right]\right\} ,
	\end{array}
	\end{equation}
	\begin{equation}\label{e3.24a}\begin{array}{l}
	f_Y(y)   =  (2(\gamma y^2+\alpha)^{-\frac{1}{2}}\times\\
	\hfill \hspace{0.8cm}\exp\left\{-\frac{1}{2}\left[2(\beta y^2-y+a_{00})
	- \displaystyle\frac{1}{2(\gamma y^2+\alpha)}\right]\right\} .
	\end{array}
	\end{equation}

\begin{figure}[ht!]
	\centering
	\caption{ Marginal density of $X$ plots for dependence, near independence and bimodality } 
	\subfloat[$\alpha=1, \beta=4, \gamma=5$ ]{\includegraphics[width=0.35\columnwidth]{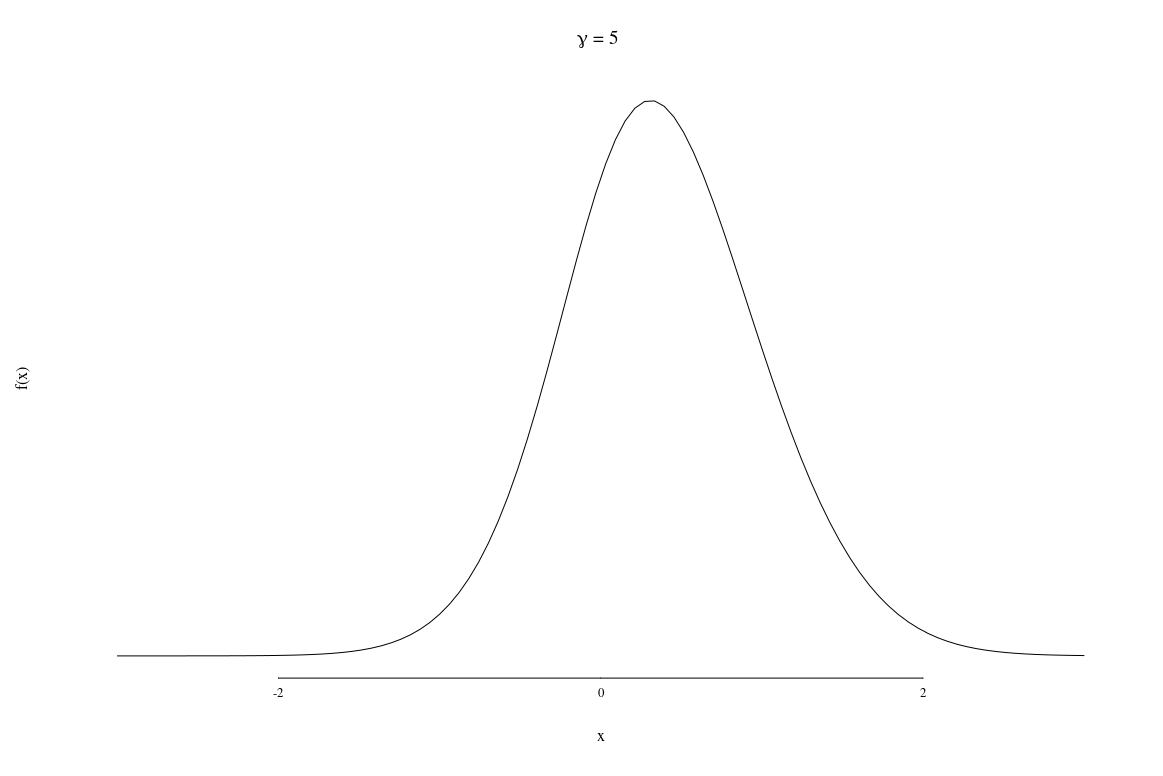}} 
	\subfloat[$\alpha=1, \beta=4, \gamma=0.12$]{\includegraphics[width=0.35\columnwidth]{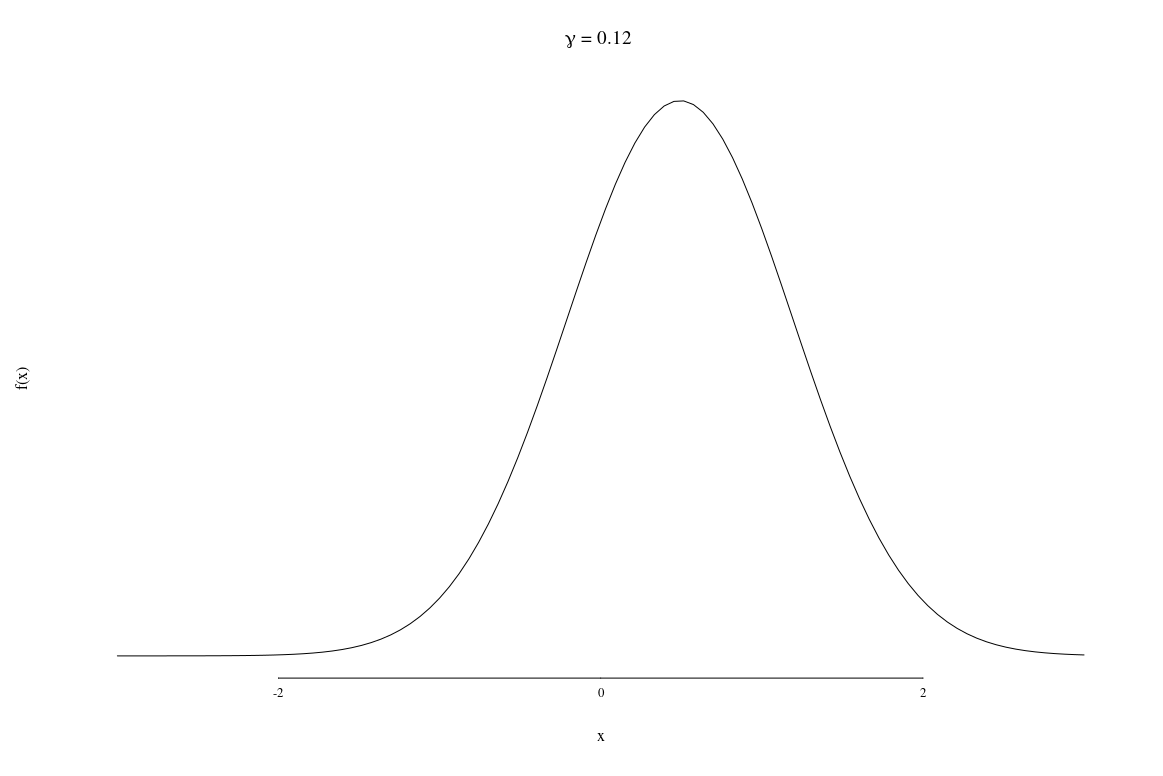}} 
	\subfloat[$\alpha=0.099, \beta=0.088, \gamma=0.12$]{\includegraphics[width=0.35\columnwidth]{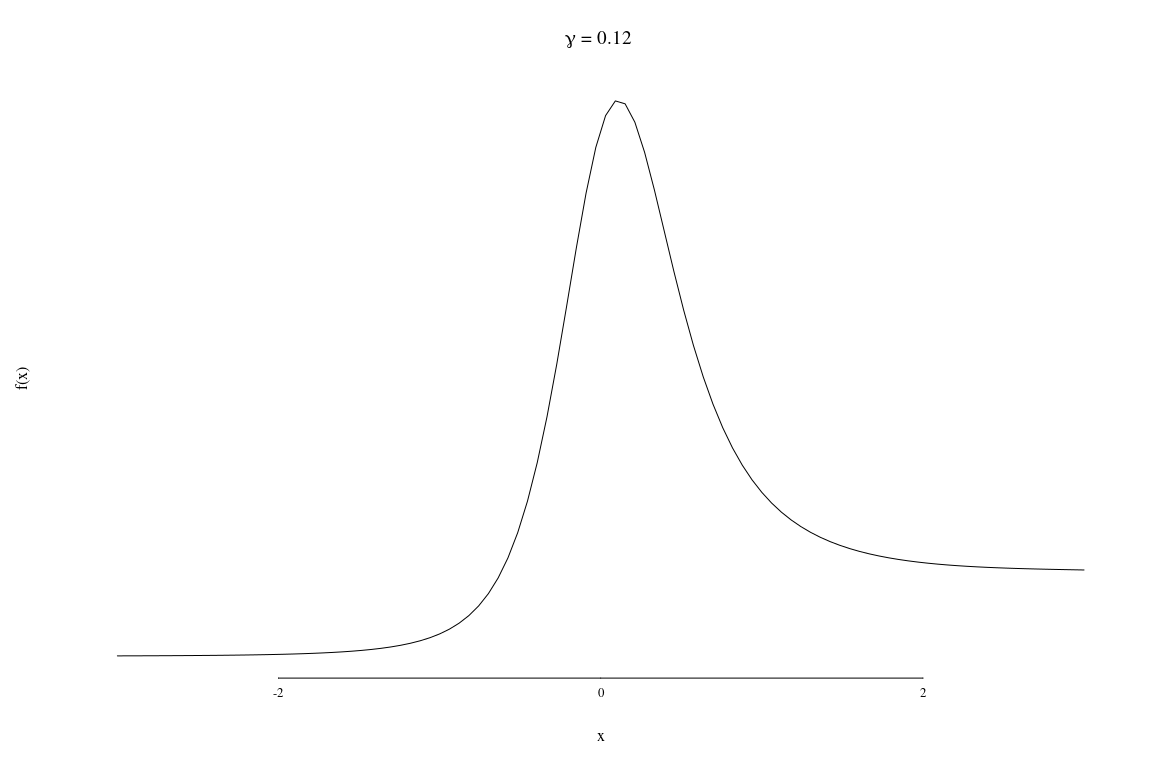}} 
	\label{marginalx}
\end{figure}

\begin{figure}[ht!]
	\centering
	\caption{ Marginal density of $Y$ plots for dependence, near independence and bimodality } 
	\subfloat[$\alpha=1, \beta=4, \gamma=5$ ]{\includegraphics[width=0.35\columnwidth]{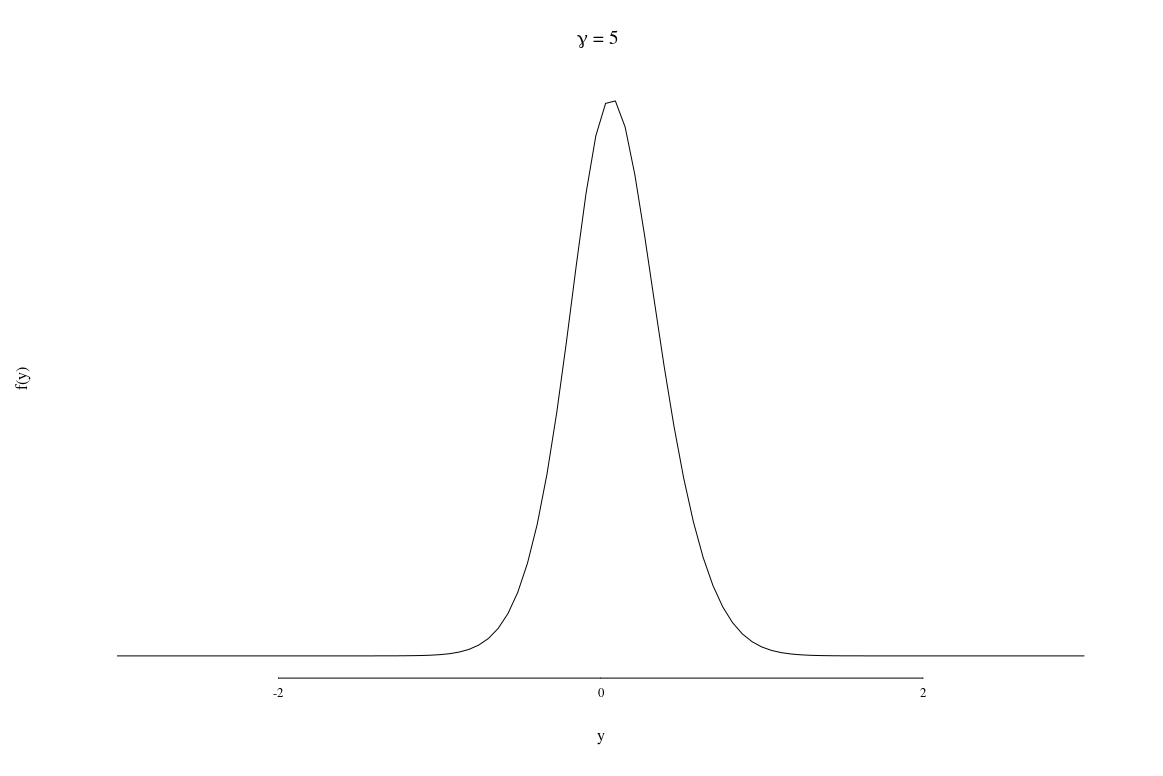}} 
	\subfloat[$\alpha=1, \beta=4, \gamma=0.12$]{\includegraphics[width=0.35\columnwidth]{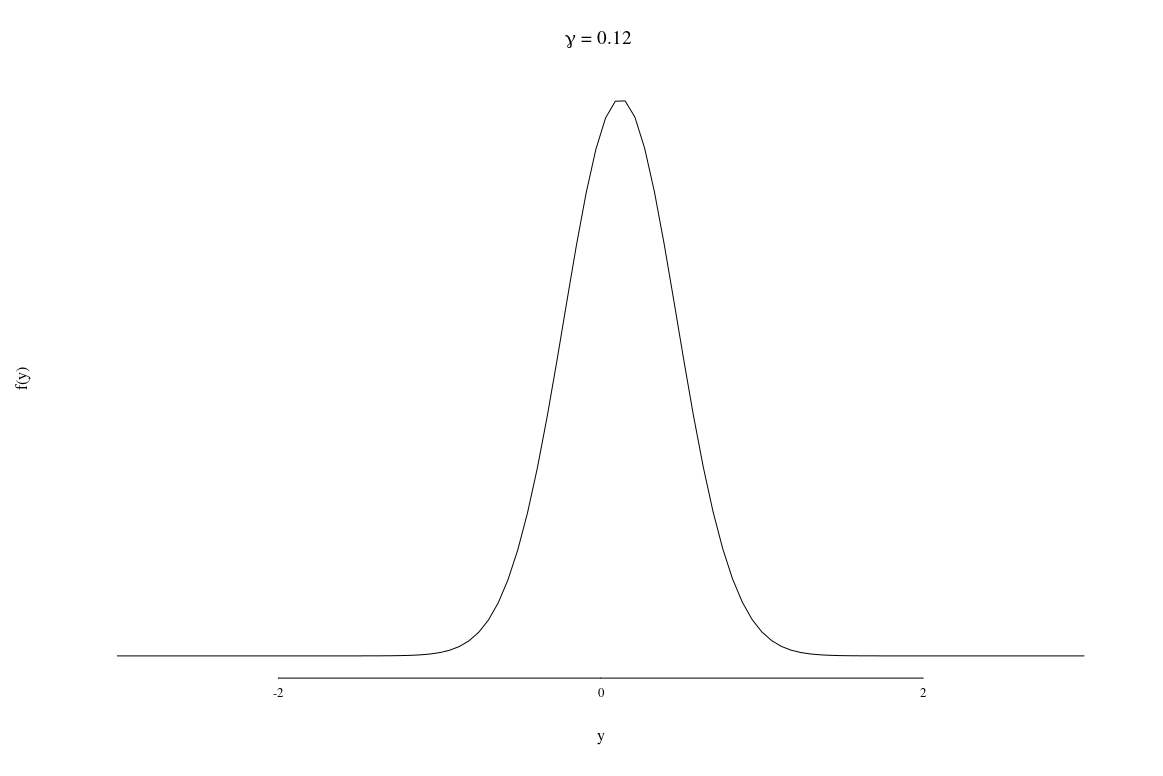}} 
	\subfloat[$\alpha=0.099, \beta=0.088, \gamma=0.12$]{\includegraphics[width=0.35\columnwidth]{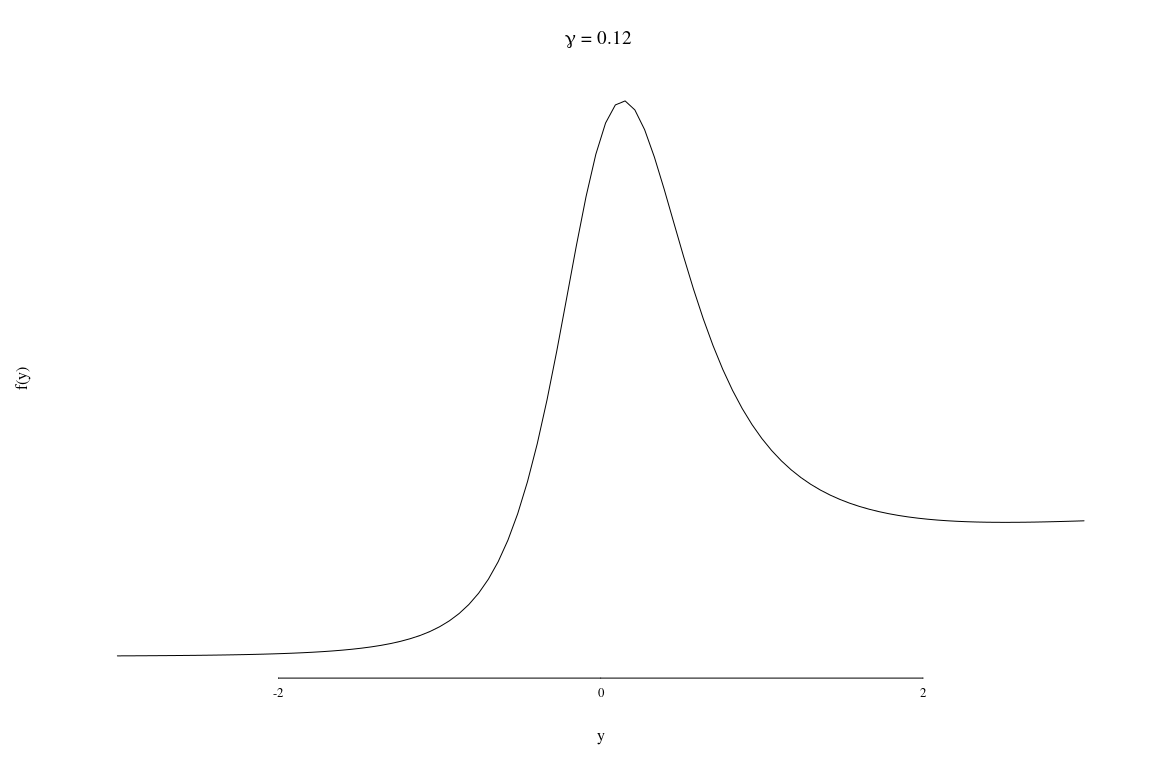}}
	\label{marginaly}
\end{figure}

	Observe that the marginal density of $X$ is a product of a normal density with mean $-1/2\alpha \neq 0$ and a function that is symmetric about $0$. This density will thus be asymmetric unless the second factor is constant, which occurs only if $\gamma=0$, i.e,. only in the case in which $X$ and $Y$ are independent. 
	Analogously, the density of $Y$ will be asymmetric except in the case of independence.   See Figure  \ref{marginalx} and Figure \ref{marginaly} for the marginal densities of $X$ and $Y$ for different choices of parameters, with strong dependence, near independence and bimodality, respectively.
	
\section{Estimation and inference}

\subsection{Maximum likelhood estimation for the bivariate equi-dispersed normal conditionals distribution}
	Consider the three parameter equi-dispersed normal density given in $(2.17)$
	\begin{equation}\label{eq-di-no-co}
		f_{X,Y}(x,y:\alpha,\beta,\gamma)\propto exp\{-[\alpha x^2+\beta y^2 +\gamma x^2y^2-x-y]\}
	\end{equation}
     where $\alpha>0,\beta>0$ and $\gamma \geq 0$. For computing maximum likelihood estimators (m.l.e) one needs to consider the complete density which involves a normalizing factor to ensure that the density integrates to $1$. To this end, we define 
         \begin{eqnarray}\label{eq-di-no-co}
     \kappa^{-1}(\alpha,\beta,\gamma)&=& \int_{-\infty}^{\infty}\int_{-\infty}^{\infty} exp\{-[\alpha x^2+\beta y^2 +\gamma x^2y^2-x-y]\} dx dy.
     \end{eqnarray}
    The bivariate equi-dispersed normal density will then be
        \begin{equation}\label{eq-di-no-co}
    	f_{X,Y}(x,y:\alpha,\beta,\gamma) = \kappa(\alpha,\beta,\gamma) \exp\{-[\alpha x^2+\beta y^2 +\gamma x^2y^2-x-y]\}, \mbox{   } -\infty <x,y<\infty.
    \end{equation}
    For the given bivariate random sample of size $n$ from the above density, i.e., $(X_1,Y_1),...,(X_n,Y_n)$, the likelihood function is
        \begin{equation}\label{eq-di-no-co-n}
    	\log(L(\alpha,\beta,\gamma)) = n \log(\kappa(\alpha,\beta,\gamma))  - \alpha 
    	\sum_{i=1}^{n} X_i^2 - \beta \sum_{i=1}^{n}Y_i^2 - \gamma \sum_{i=1}^{n}\sum_{j=1}^{n}X_i^2Y_j^2  + \sum_{i=1}^{n} X_i+\sum_{i=1}^{n}  Y_i. 
    \end{equation}
   Note the explicit expression for the maximum likelihood estimators are not possible and one needs to depend on a numerical method to find the maximum likelihood estimates for the given data. 
   
   \begin{remark}
   	We make a remark that the algorithm needs to evaluate the normalizing factor $\kappa(.)$,  for each choice of parameter values. For the R software this optimization can be handled by defining "closure" and the corresponding code is included in Appendix C.  Also, due to the special optimization process we recommend the "Rvmmin" algorithm for maximization. We refer to Nash \cite{n19} for further details on numerical optimization for nested functions.
   \end{remark}

	\subsection{On pseudo-likelihood estimation for the bivariate \\ equi-dispersed normal conditionals distribution}
	Instead of using the likelihood function, which, as we have seen, is challenging to maximize, it is natural consider pseudo-likelihood as a more convenient alternative to obtain consistent estimates which in general are somewhat less efficient than maximum likelihood estimates, were they available. A convenient introduction to pseudo-likelihood estimation may be found in Arnold and Strauss (1991). The pseudo likeihood function corresponding to a bivariate sample $(X_i,Y_i), \ \ i=1,2,...,n$ from the density 
	$f_{X,Y}(x,y;\underline{\theta})$ is given by
	$$PL(\underline{\theta})=\prod_{i=1}^n f_{X|Y}(x_i|y_i;\underline{\theta})f_{Y|X}(y_i|x_i;\underline{\theta}).$$
For the case of distributions with equi-dispersed normal conditionals, the log-pseudo-likelihood is given by.
\begin{eqnarray}\label{log-pseu-lik}
log PL(\alpha,\beta,\gamma)=c+(1/2)\sum_{i=1}^n log(\gamma y_i^2 +\alpha)-\sum_{i=1}^n\frac{1}{4(\gamma y_i^2+\alpha)}-\sum_{i=1}^n x_i^2(\gamma y_i^2+\alpha)  \nonumber \\
+(1/2)\sum_{i=1}^n log(\gamma x_i^2 +\beta)-\sum_{i=1}^n\frac{1}{4(\gamma x_i^2+\beta)}-\sum_{i=1}^n y_i^2(\gamma x_i^2+\beta).
\end{eqnarray}
The pseudo-likelihood estimates of $\alpha,\beta$ and $\gamma$ are the values of these parameters that maximize the pseudo-likelihood which can be achieved by considering the log-pseudo-likelihood
 in (\ref{log-pseu-lik}) above, 

Differentiating with respect to $\alpha, \beta$ and $\gamma$ and equating to $0$ yields the following pseudo-likelihood equations.

\begin {eqnarray}
&&(1/2)\sum_{i=1}^n 1/(\gamma y_i^2 +\alpha)+\sum_{i=1}^n\frac{1}{4(\gamma y_i^2+\alpha)^2}=\sum_{i=1}^n x_i^2,  \label{PL1} \\
&&(1/2)\sum_{i=1}^n 1/(\gamma x_i^2 +\beta)+\sum_{i=1}^n\frac{1}{4(\gamma x_i^2+\beta)^2}=\sum_{i=1}^n y_i^2 ,   \label{PL2} \\
&&(1/2)\sum_{i=1}^n y_i^2/(\gamma y_i^2 +\alpha)+\sum_{i=1}^n\frac{y_i^2}{4(\gamma y_i^2+\alpha)^2}  \nonumber \\
&&+(1/2)\sum_{i=1}^n x_i^2 /(\gamma x_i^2 +\beta)+\sum_{i=1}^n\frac{x_i^2}{4(\gamma x_i^2+\beta)^2}=2\sum_{i=1}^n x_i^2 y_i^2 . \label{PL3}
\end{eqnarray}

Note that the left hand sides of these equations are well-behaved. For a fixed value of $\gamma$  the left side of equation (\ref{PL1}) is a decreasing function of $\alpha$. For a fixed value of $\gamma$  the left side of equation (\ref{PL2}) is a decreasing function of $\beta$, and for  fixed value of $\alpha$ and $\beta$  the left side of equation (\ref{PL3}) is a decreasing function of $\gamma$. As a consequence an iterative scheme can be used to identify the corresponding pseudo-likeliood estimates.
	
	\subsection{Likelihood ratio test for bivariate equi-dispersed normal conditionals}
	It is natural to consider an ordinary bivariate normal distribtion as a $5$ parameter alternative to the $3$ parameter equidispered conditionals model, both of which are nested within the $8$-parameter normal conditionals model with density (\ref{e3.26}). There is very little overlap between the classical normal model and the model with equi-dispersed normal conditionals. The only distributions that are  in both families are those with independent equi-dispersed normal marginals. It is possible to envision a likelihood ratio test for the equi-dispersed normal conditionals model within the full $8$-parameter normal conditionals model, but the effort will require non-trivial computer intensive maximum likelihood estimation of the parameters in the models. It will of course be possible to compare the various models using an AIC or BIC criterion, perhaps  using pseudo-likelihood parameter estimates for the $3$-parameter model.

		\section{Examples, simulated and real-world}
In the following three sub-sections we provide a bootstrapped simulation study of the m.l.e's and pseudo m.l.e's of the parameters  of the bivariate density given in (\ref{eq-di-no-co}) and also include two examples of real-life application of the proposed model.

	\subsection{Simulated data}
 A simple simulation algorithm for the bivariate equi-dispersed normal conditionals model, for a given $\alpha$, $\beta$ and $\gamma$, involves the following steps.
	\begin{description}
		\item[Step 1:] Simulate $x$ from the marginal density given in  (\ref{X-marg}). Note that for the given parameter values the normalizing constant $a_{00}$ is fixed and is computed by numerical integration.
		\item[Step 2:] Next, for the given $x$ simulate $y$ from a $N\Big (\frac {1}{2(\gamma x^2+\beta)}, \frac {1}{2(\gamma x^2+\beta)}\Big )$ distribution.
	\end{description}
	Repeat the above two steps for the desired number of observations.
	\bigskip
	
\begin{figure}%
		\centering
		\subfloat[p.m.l.e for $\alpha(=1)$]{{\includegraphics[width=70mm]{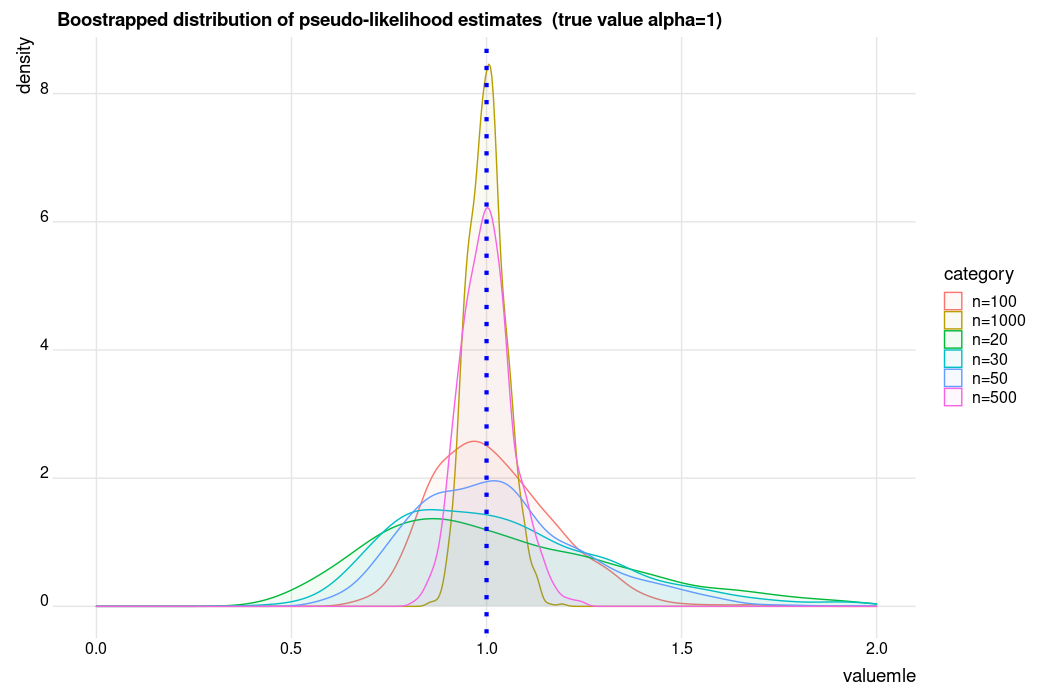} }}%
		\qquad
		\subfloat[ m.l.e for $\alpha(=1)$]{{\includegraphics[width=70mm]{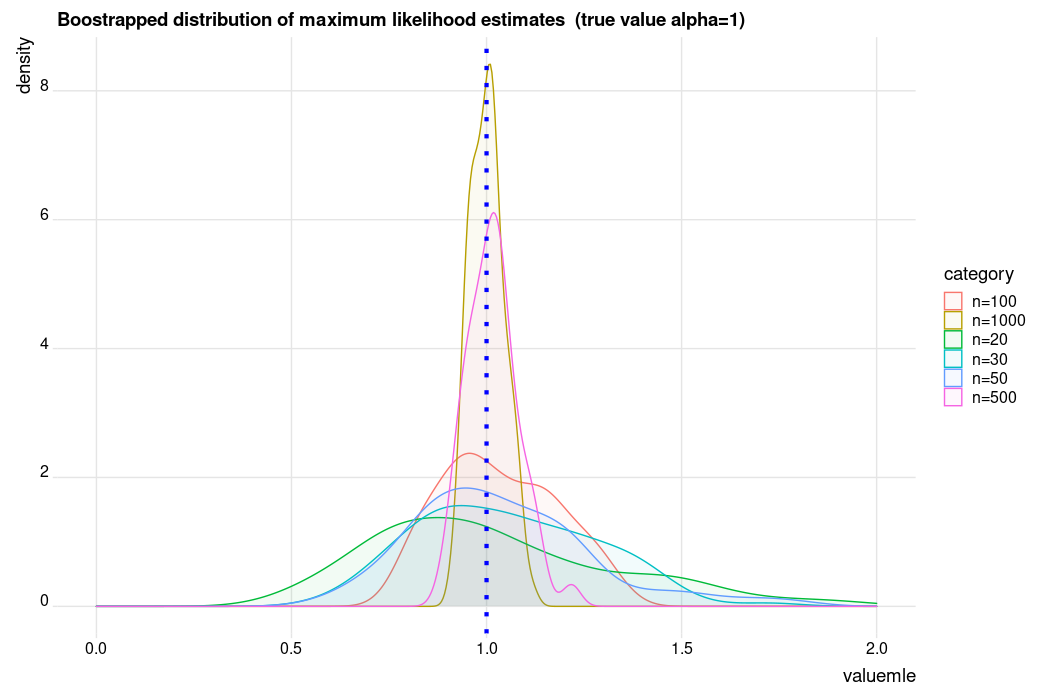} }}%
		\caption{$\alpha=1$,$\beta=4$, $\gamma=5$ (dependence)}%
		\label{fig1}%
	\end{figure}

\begin{figure}%
		\centering
		\subfloat[p.m.l.e for $\beta(=4)$]{{\includegraphics[width=70mm]{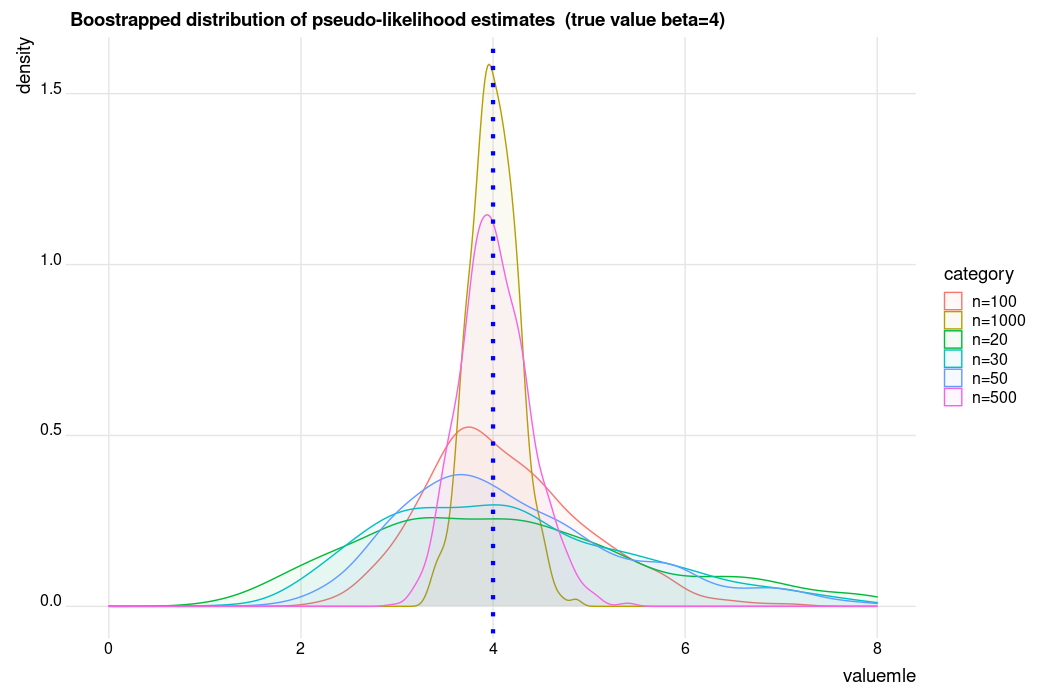} }}%
		\qquad
		\subfloat[ m.l.e for $\beta(=4)$]{{\includegraphics[width=70mm]{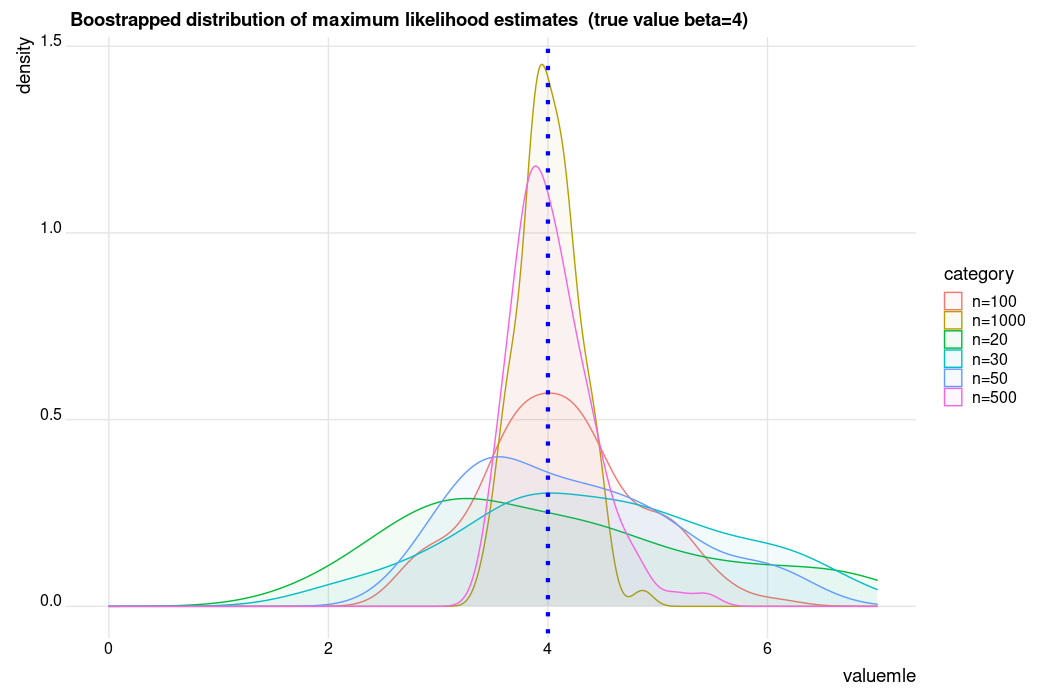} }}%
		\caption{$\alpha=1$,$\beta=4$, $\gamma=5$(dependence)}%
		\label{fig2}%
	\end{figure}

\begin{figure}%
		\centering
		\subfloat[p.m.l.e for $\gamma(=5)$]{{\includegraphics[width=70mm]{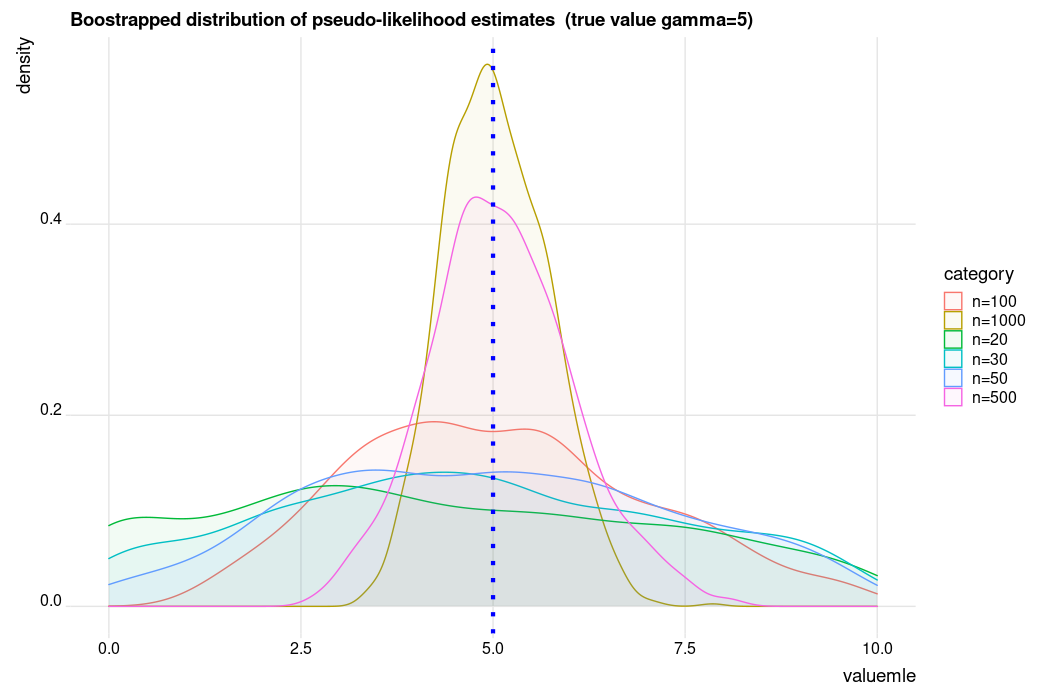} }}%
		\qquad
		\subfloat[ m.l.e for $\gamma(=5)$]{{\includegraphics[width=70mm]{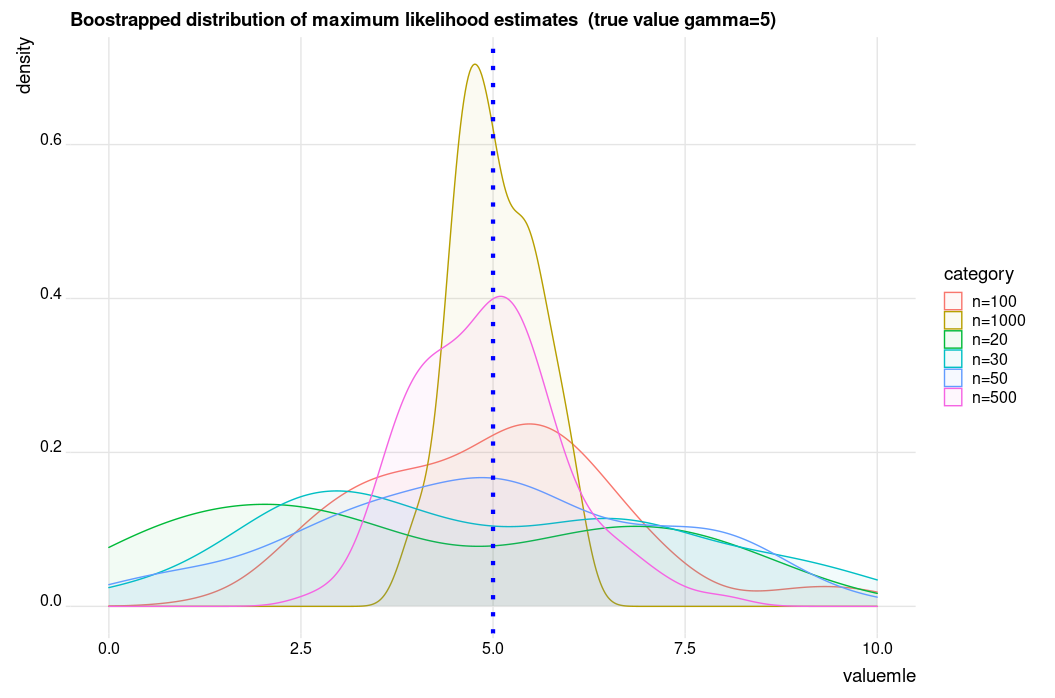} }}%
		\caption{$\alpha=1$,$\beta=4$, $\gamma=5$(dependence)}%
		\label{fig3}%
	\end{figure}

\begin{figure}%
	\centering
	\subfloat[p.m.l.e for $\alpha(=1)$]{{\includegraphics[width=70mm]{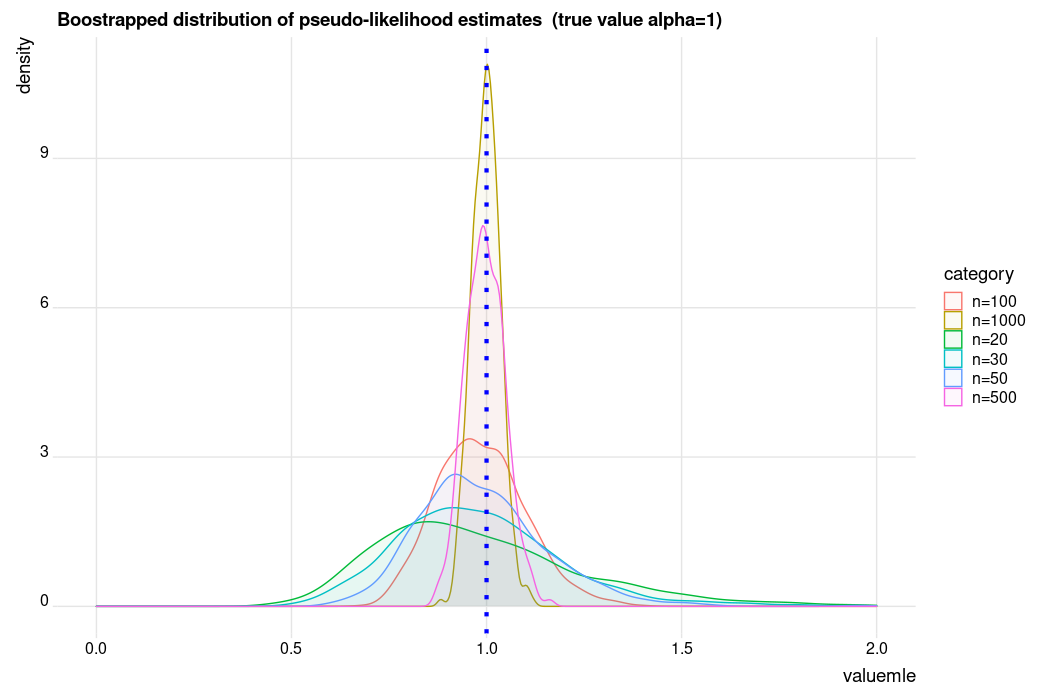} }}%
	\qquad
	\subfloat[ m.l.e for $\alpha(=1)$]{{\includegraphics[width=70mm]{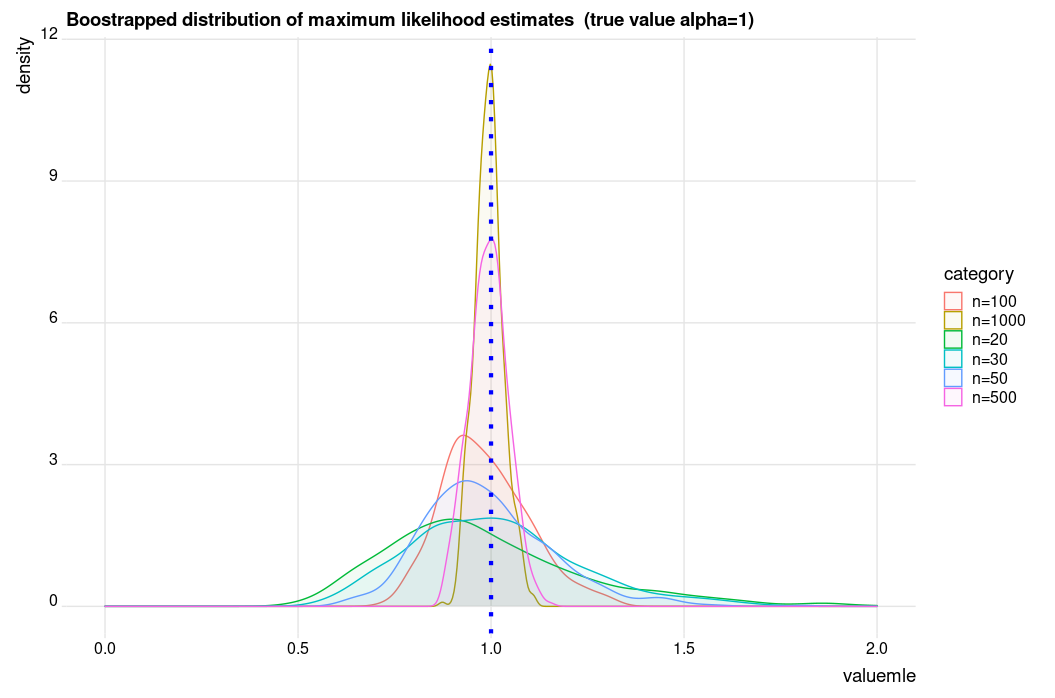} }}%
	\caption{$\alpha=1$,$\beta=4$, $\gamma=0.12$(close to independence)}%
	\label{fig4}%
\end{figure}

\begin{figure}%
	\centering
	\subfloat[p.m.l.e for $\beta(=4)$]{{\includegraphics[width=70mm]{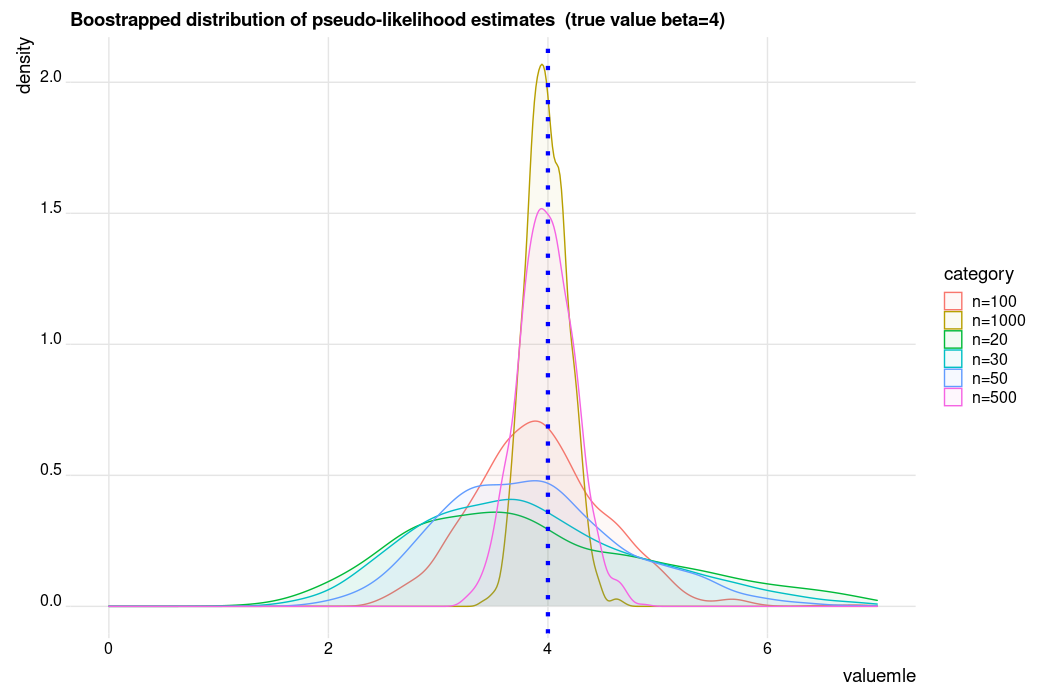} }}%
	\qquad
	\subfloat[ m.l.e for $\beta(=4)$]{{\includegraphics[width=70mm]{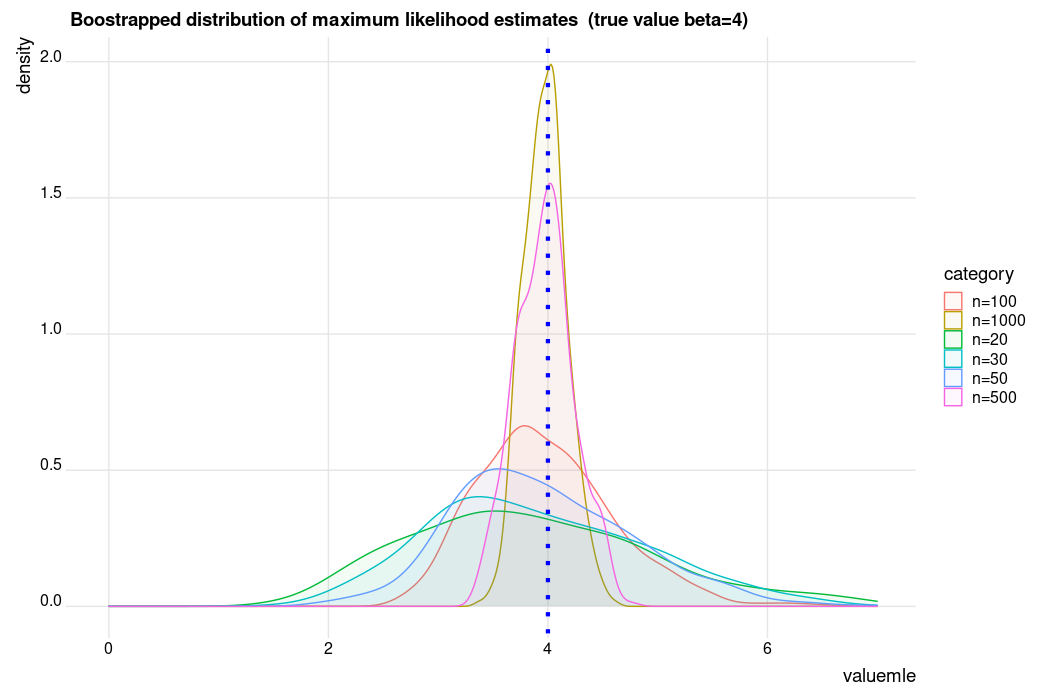} }}%
	\caption{$\alpha=1$,$\beta=4$, $\gamma=0.12$(close to independence)}% 
	\label{fig5}%
\end{figure}

\begin{figure}%
	\centering
	\subfloat[p.m.l.e for $\gamma(=0.12)$]{{\includegraphics[width=70mm]{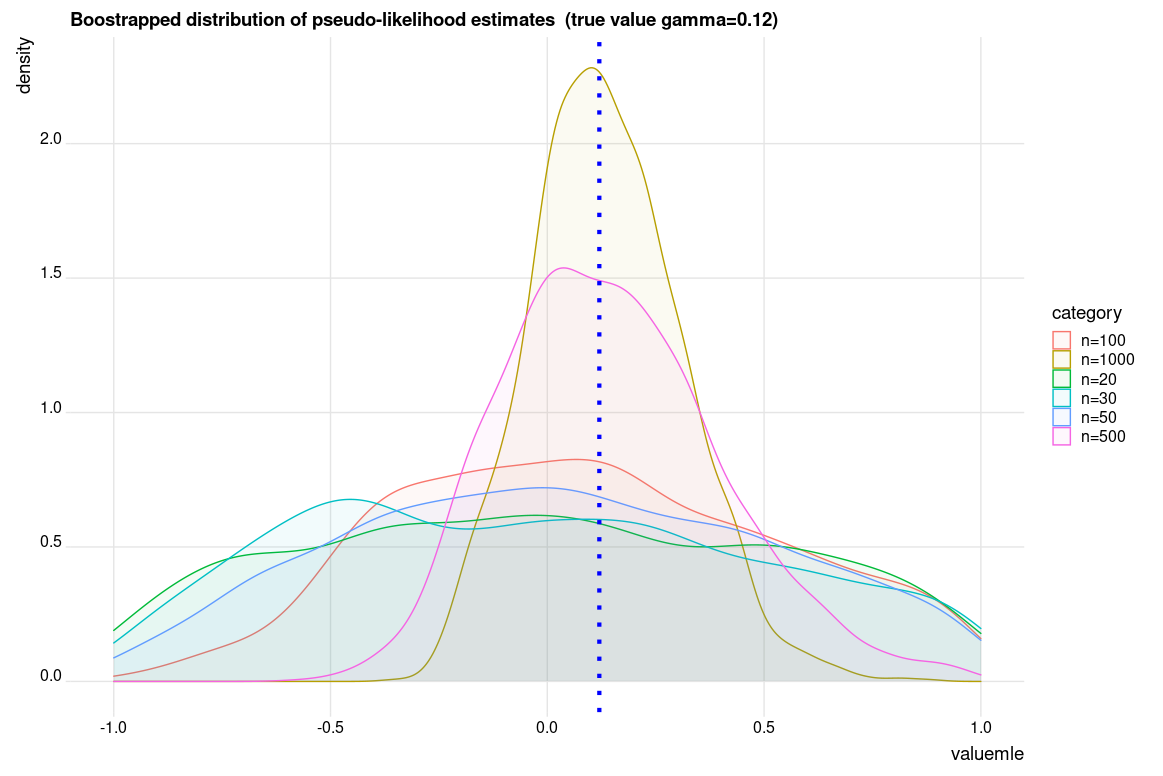} }}%
	\qquad
	\subfloat[ m.l.e for $\gamma(=0.12)$]{{\includegraphics[width=70mm]{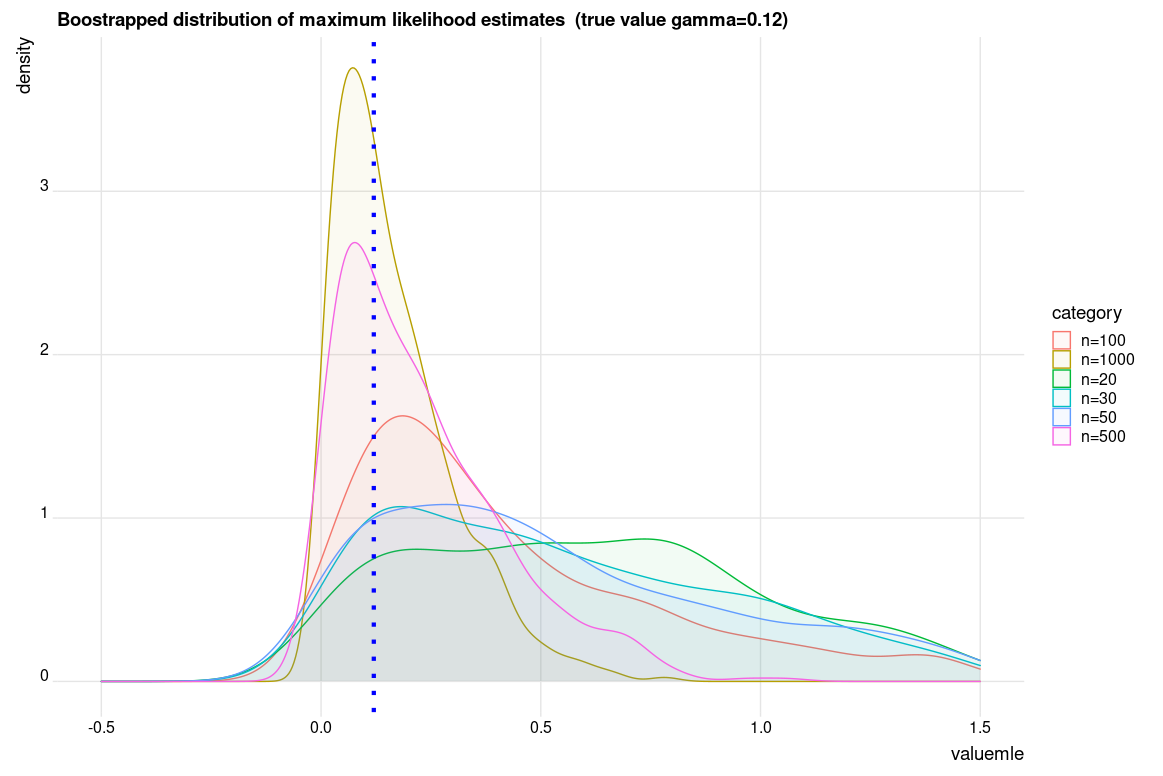} }}%
	\caption{$\alpha=1$,$\beta=4$, $\gamma=0.12$(close to independence)}%
	\label{fig6}%
\end{figure}

	We have simulated  $5000$ data sets of sample size $n= 20,30,50,100,$
$500,1000$  from the density in  (\ref{eq-di-no-co}) for two different parametric configurations.  We refer to Figures \ref{fig1}--\ref{fig6}  for the boostrapped distribution of the
pseudo m.l.e and m.l.e's. and  also, see Tables \ref{sg5} \& \ref{sg012} for summary values from the boostrapped samples (includes mean, standard error(SD)  and $95\%$ confidence intervals).

\begin{table}
	\caption{Simulation study:  true values are $\alpha=1$, $\beta =4$ and $\gamma =5$ }  
	\label{sg5}
	\small % text size of table content
	\centering % center the table
	\begin{tabular}{lcccccccr} % alignment of each column data
		\toprule[\heavyrulewidth]\toprule[\heavyrulewidth]
		\textbf{$n$ }  & \textbf{P} & \textbf{MLE} &  \textbf{SE(MLE)} &   \textbf{95\% CI (MLE)} & \textbf{PMLE} &  \textbf{SE(PMLE)} &   \textbf{95\% CI (PMLE)} \\ 
		\midrule
		\multirow{3}{*}{$20$} & $\alpha$ & $1.042$ & $0.351$ & $(0.559,1.885)$ & $1.105$ & $0.434$ & $(0.562,2.163)$  \\
		& $\beta$ &  $4.603$ & $2.095$ & $(2.103,9.188)$   & $4.546$ & $2.352$ & $(1.770,10.250)$  \\
		& $\gamma$ & $5.601$ & $4.609$ & $(0.000,14.494)$   & $7.165$ & $6.886$ & $(0.000,24.870)$  \\
	
		\bottomrule[\heavyrulewidth] 
			\multirow{3}{*}{$30$} & $\alpha$ & $1.084$ & $0.292$ & $(0.686,1.923)$ & $1.060$ & $0.298$ & $(0.616,1.777)$\\
		& $\beta$ &  $4.501$ & $1.257$ & $(2.141,6.809)$   & $4.342$ & $1.736$ & $(2.099,8.520)$  \\
		& $\gamma$ & $5.986$ & $4.134$ & $(0.753,15.677)$   & $6.473$ & $4.990$ & $(0.075,18.209)$  \\
		
		\bottomrule[\heavyrulewidth] 
			\multirow{3}{*}{$50$} & $\alpha$ & $1.051$ & $0.265$ & $(0.673,1.690)$ & $1.034$ & $0.223$ & $(0.677,1.548)$ \\
		& $\beta$ &  $4.355$ & $1.227$ & $(2.837,7.495)$   & $4.230$ & $1.253$ & $(2.337,7.174)$  \\
		& $\gamma$ & $5.386$ & $2.768$ & $(0.688,11.324)$   & $5.776$ & $3.545$ & $(0.441,13.547)$  \\
		
		\bottomrule[\heavyrulewidth] 
			\multirow{3}{*}{$100$} & $\alpha$ & $1.034$ & $0.147$ & $(0.789,1.304)$ & $1.010$ & $0.152$ & $(0.758,1.354)$  \\
		& $\beta$ &  $4.114$ & $0.700$ & $(2.823,5.498)$   & $4.106$ & $0.806$ & $(2.772,5.931)$  \\
		& $\gamma$ & $5.104$ & $1.822$ & $(2.394,9.561)$   & $5.309$ & $2.154$ & $(1.782,9.736)$  \\
		
		\bottomrule[\heavyrulewidth] 
			\multirow{3}{*}{$500$} & $\alpha$ &  $1.015$ & $0.066$ & $(0.911,1.140)$ & $1.006$ & $0.182$ & $(0.881,1.137)$ \\
		& $\beta$ &  $4.047$ & $0.368$ & $(3.519,4.861)$   & $4.018$ & $0.347$ & $(3.392,4.744)$  \\
		& $\gamma$ & $4.936$ & $0.971$ & $(3.366,6.994)$   & $5.105$ & $0.950$ & $(3.321,7.235)$  \\
		
		\bottomrule[\heavyrulewidth] 
			\multirow{3}{*}{$1000$} & $\alpha$ & $1.002$ & $0.044$ & $(0.932,1.084)$ & $1.001$ & $0.048$ & $(0.914,1.094)$\\
		& $\beta$ &  $4.009$ & $0.268$ & $(3.514,4.496)$   & $4.003$ & $0.250$ & $(3.538,4.519)$  \\
		& $\gamma$ &  $5.040$ & $0.536$ & $(4.048,6.064)$   & $5.045$ & $0.671$ & $(3.900,6.387)$  \\
		
		\bottomrule[\heavyrulewidth] 
		\bigskip

	\end{tabular}
\end{table}

\begin{table}
	\caption{Simulation study:  true values are $\alpha=1$, $\beta =4$ and $\gamma =0.12$ }  
	\label{sg012}
	\small % text size of table content
	\centering % center the table
	\begin{tabular}{lcccccccr} % alignment of each column data
		\toprule[\heavyrulewidth]\toprule[\heavyrulewidth]
		\textbf{$n$ }  & \textbf{P} & \textbf{MLE} &  \textbf{SE(MLE)} &   \textbf{95\% CI (MLE)} & \textbf{PMLE} &  \textbf{SE(PMLE)} &   \textbf{95\% CI (PMLE)} \\ 
		\midrule
		\multirow{3}{*}{$20$} & $\alpha$ & $0.987$ & $0.266$ & $(0.617,1.641)$ & $1.058$ & $0.330$ & $(0.607,1.772)$  \\
		& $\beta$ &  $4.041$ & $1.451$ & $(2.064,7.836)$   & $4.342$ & $1.724$ & $(2.037,9.162)$  \\
		& $\gamma$ & $1.043$ & $1.821$ & $(0.000,6.585)$   & $0.871$ & $2.586$ & $(0.000,6.594)$  \\
		
		\bottomrule[\heavyrulewidth] 
		\multirow{3}{*}{$30$} & $\alpha$ & $1.006$ & $0.213$ & $(0.649,1.514)$ & $1.034$ & $0.264$ & $(0.646,1.646)$\\
		& $\beta$ &  $3.954$ & $1.069$ & $(2.259,6.332)$   & $4.266$ & $1.408$ & $(2.266,6.946)$  \\
		& $\gamma$ & $0.813$ & $1.247$ & $(0.000,3.938)$   & $0.568$ & $1.529$ & $(0.000,4.286)$  \\
		
		\bottomrule[\heavyrulewidth] 
		\multirow{3}{*}{$50$} & $\alpha$ & $0.990$ & $0.165$ & $(0.731,1.395)$ & $1.007$ & $0.184$ & $(0.692,1.409)$ \\
		& $\beta$ &  $3.967$ & $0.835$ & $(2.581,5.734)$   & $4.066$ & $1.016$ & $(2.582,6.059)$  \\
		& $\gamma$ & $0.607$ & $0.903$ & $(0.000,3.331)$   & $0.460$ & $1.080$ & $(0.000,3.377)$  \\
		
		\bottomrule[\heavyrulewidth] 
		\multirow{3}{*}{$100$} & $\alpha$ & $0.986$ & $0.112$ & $(0.793,1.238)$ & $1.001$ & $0.128$ & $(0.785,1.276)$  \\
		& $\beta$ &  $3.980$ & $0.619$ & $(2.982,5.321)$   & $4.043$ & $0.653$ & $(2.886,5.402)$  \\
		& $\gamma$ & $0.393$ & $0.551$ & $(0.000,1.847)$   & $0.314$ & $0.667$ & $(0.000,1.962)$  \\
		
		\bottomrule[\heavyrulewidth] 
		\multirow{3}{*}{$500$} & $\alpha$ &  $0.994$ & $0.050$ & $(0.900,1.091)$ & $1.001$ & $0.109$ & $(0.898,1.113)$ \\
		& $\beta$ &  $3.971$ & $0.264$ & $(3.484,4.515)$   & $4.002$ & $0.413$ & $(3.466,4.590)$  \\
		& $\gamma$ & $0.198$ & $0.221$ & $(0.000,0.674)$   & $0.268$ & $0.435$ & $(0.000,0.723)$  \\
		
		\bottomrule[\heavyrulewidth] 
		\multirow{3}{*}{$1000$} & $\alpha$ & $0.995$ & $0.037$ & $(0.927,1.072)$ & $1.000$ & $0.036$ & $(0.929,1.077)$\\
		& $\beta$ &  $3.975$ & $0.197$ & $(3.613,4.371)$   & $4.000$ & $0.200$ & $(3.593,4.419)$  \\
		& $\gamma$ &  $0.166$ & $0.151$ & $(0.000,0.525)$   & $0.155$ & $0.155$ & $(0.000,0.529)$  \\
		
		\bottomrule[\heavyrulewidth] 
		\bigskip

	\end{tabular}
\end{table}

\bigskip

\bigskip

We summarize Tables  \ref{sg5} and \ref{sg012} by the following general remarks.  We note that with an increase in sample size, both the pseudo and the actual m.l.e's standard errors (SE) decrease. Also, the $95\%$ confidence intervals using the actual m.l.e's have shorter length compared to the confidence intervals constructed using the pseudo m.l.e's.  In particular, we observe that,  for the sample size greater than $30$ the m.l.e's  approach the true parameter values with decreasing standard errors.  The corresponding pseudo m.l.e's behave in a similar fashion as sample sizes increase but have higher standard errors.   We also make a remark   that for values of  $\gamma$ close to zero both  pseudo and actual m.l.e's algorithms fail to converge in many cases for small sample sizes.  Finally, we also recommend that the pseudo m.l.e's can be considered as the primary choice for the initial values for the numerical computation of actual  m.l.e's for any sample size.

	\subsection{Real-life data I}
	In the following we considered Piedmont wines data on chemical properties of 178 specimens of three types of wine produced in the Piedmont region of Italy. The data represent 27 chemical measurements on each of 178 wine specimens belonging to three types of wine produced in the Piedmont region of Italy. The measurements on three types of wines, includes, alcohol (alcohol percentage),sugar	(sugar-free extract), uronic (uronic acids), hue (numerical), nitrogen (total nitrogen), methanol (methanol),  etc.
	 We refer to Forina et al. \cite{f86} and Azzalini \cite{az22} for further reference on the  Piedmont data.

    Here we consider two measurements,  Uronic acids ($X$) and Hue ($Y$)  on the three types of wine produced in the Piedmont, see Figure 15 for a scatter plot of this bivariate data set.

\begin{figure}[ht!]
	\centering
	\includegraphics[width=150mm]{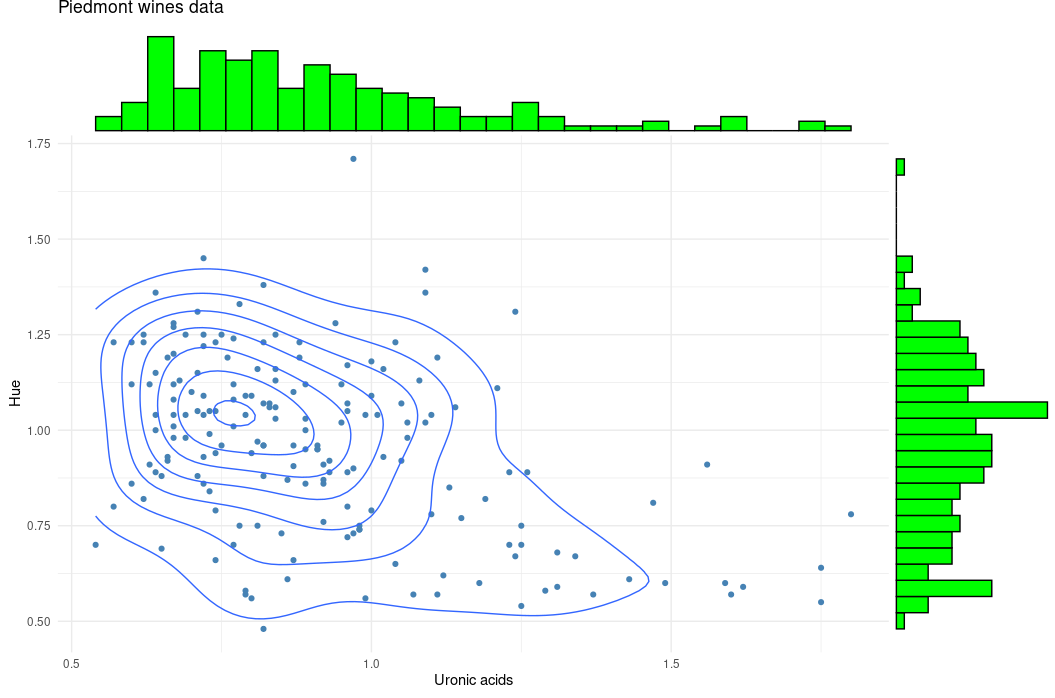}	
	\label{wuh}
	\caption{Piedmont wines data: Uronic acids ($X$) and Hue ($Y$) scatter plot}
\end{figure}

In the following we fit four models for the above bivariate data:
\begin{itemize}
 
	\item Model I (dependent):  Here we considered the dependent equi-disperesed conditionals normal model. We refer to  Table \ref{fullD} for the  fitted m.l.e's (both pseudo and actual) and AIC, respectively.
	\item Model II (indepedent): Here we considered the independent equi-disperesed normal model. We refer to  Table \ref{fullI} for the  fitted m.l.e's (both pseudo and actual) and AIC, respectively.
	\item Model III (bivariate normal): Here we considered the classical bivariate normal model. We refer to  Table \ref{fullBD} for the  fitted m.l.e's  and AIC, respectively.
	\item Model IV (bivariate normal indepedent):Here we considered the bivariate normal model with indepedent marginals. We refer to  Table \ref{fullBI} for the  fitted m.l.e's  and AIC, respectively.
\end{itemize}

	\begin{table}[ht!]
	\caption{Model I on  Uronic acids ($X$) and Hue ($Y$) }  
	\label{fullD}
	\small % text size of table content
	\centering % center the table
	\begin{tabular}{lcccccr} % alignment of each column data
		\toprule[\heavyrulewidth]\toprule[\heavyrulewidth]
		\textbf{$n$ }  & \textbf{P} & \textbf{MLE} & \textbf{PMLE} &  \textbf{AIC } \\ 
		\midrule
		\multirow{3}{*}{$178$} & $\alpha$ & $0.829$ & $0.575$ & \multirow{3}{*}{$-566.656$} \\
		& $\beta$ &  $0.786$ & $0.559$ \\
		& $\gamma$ & $0.051$ & $0.338$\\
		
		\bottomrule[\heavyrulewidth] 
		\bigskip
		
	\end{tabular}
\end{table}

	\begin{table}[ht!]
	\caption{Model II on  Uronic acids ($X$) and Hue ($Y$) }  
	\label{fullI}
	\small % text size of table content
	\centering % center the table
	\begin{tabular}{lcccccr} % alignment of each column data
		\toprule[\heavyrulewidth]\toprule[\heavyrulewidth]
		\textbf{$n$ }  & \textbf{P} & \textbf{MLE} & \textbf{PMLE} &  \textbf{AIC } \\ 
		\midrule
		\multirow{2}{*}{$178$} & $\alpha$ & $0.874$ & $0.828$ & \multirow{2}{*}{$-569.492$} \\
		& $\beta$ &  $0.874$ & $0.828$ \\
		
		\bottomrule[\heavyrulewidth] 
		\bigskip
		
	\end{tabular}
\end{table}

\begin{table}[ht!]
	\caption{Model III on Uronic acids ($X$) and Hue ($Y$) }  
	\label{fullBD}
	\small % text size of table content
	\centering % center the table
	\begin{tabular}{lcccccr} % alignment of each column data
		\toprule[\heavyrulewidth]\toprule[\heavyrulewidth]
		\textbf{$n$ }  & \textbf{P} & \textbf{MLE} &  \textbf{AIC } \\ 
		\midrule
		\multirow{5}{*}{$178$} & $\mu_1$ &  $0.915$ & \multirow{5}{*}{$55.505$} \\
		& $\mu_2$ & $0.957$ \\
		& $\sigma^2_1$  & $0.063$\\
			& $\sigma^2_2$  & $0.052$\\
				& $cov$ &   $-0.025$\\
		
		\bottomrule[\heavyrulewidth] 
		\bigskip
		
	\end{tabular}
\end{table}

 \begin{table}[ht!]
 	\caption{Model IV on Uronic acids ($X$) and Hue ($Y$) }  
 	\label{fullBI}
 	\small % text size of table content
 	\centering % center the table
 	\begin{tabular}{lcccccr} % alignment of each column data
 		\toprule[\heavyrulewidth]\toprule[\heavyrulewidth]
 		\textbf{$n$ }  & \textbf{P} & \textbf{MLE} &  \textbf{AIC } \\ 
 		\midrule
 		\multirow{4}{*}{$178$} & $\mu_1$ &  $0.915$ & \multirow{4}{*}{$13.563$} \\
 		& $\mu_2$ & $0.957$ \\
 		& $\sigma^2_1$  & $0.058$\\
 		& $\sigma^2_2$  & $0.052$\\
 	 		
 		\bottomrule[\heavyrulewidth] 
 		\bigskip
 		
 	\end{tabular}
 \end{table}

Note that using the AIC criterion, for the Piedmont wines data with measurements on Uronic acids ($X$) and Hue ($Y$), we recommend the bivariate dependent equi-dispersed normal conditionals model.

\subsection{Real-life data II}
In the following, we considered Australian Institute of Sport data  on 102 male and 100 female athletes collected at the Australian Institute of Sport, courtesy of Richard Telford and Ross Cunningham.   The data consist of 202 observations on 13 variables, including sex, height(cm), weight (kg), body mass index, lean body mass, red cell count, weight cell count etc.  We refer Forina et al. \cite{cw94} and Azzalini \cite{az22} for further references on the Australian Institute of Sport data.

Here we consider two variables from  the Australian Institute of Sport data, i.e.,  Body Mass Index ($X$)  and  Lean Body Mass ($Y$),  see Figure 16 for the corresponding scatter plot.

\begin{figure}[ht!]
	\centering
	\includegraphics[width=150mm]{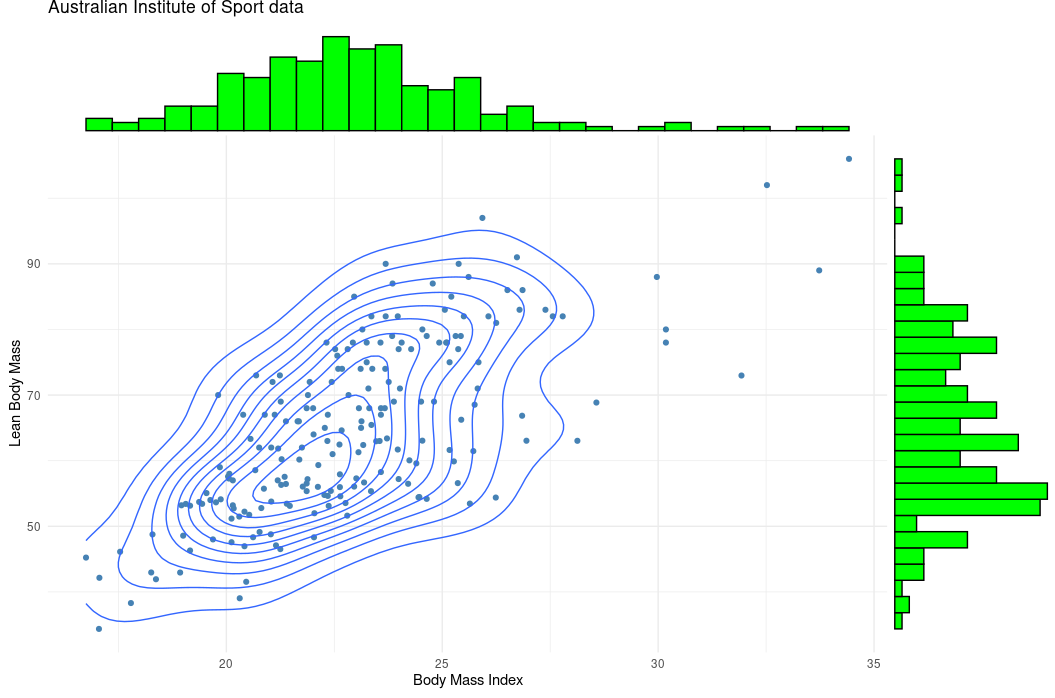}	
	\label{II}
	\caption{Australian Institute of Sport data: Body Mass Index ($X$) and Lean Body Mass ($Y$) scatter plot}
\end{figure}

\vspace{0.5cm}

For this bivariate data set we fit the following three models:
\begin{itemize}
	
	\item Model I (dependent):  Here we considered dependent equi-disperesed conditionals normal model. We refer to  Table \ref{aisD} for the  fitted m.l.e's (both pseudo and actual) and AIC, respectively.
	\item Model II (indepedent): Here we considered the model with independent equi-disperesed normal  marginals. We refer to  Table \ref{aisI} for the  fitted m.l.e's (both pseudo and actual) and AIC, respectively.
	\item Model III (bivariate normal indepedent):Here we considered the bivariate bormal model with independent marginals. We refer to  Table \ref{aisBI} for the  fitted m.l.e's  and AIC, respectively.
\end{itemize}

\begin{table}[ht!]
	\caption{Model I on  Body Mass Index ($X$) and Lean Body Mass ($Y$) }  
	\label{aisD}
	\small % text size of table content
	\centering % center the table
	\begin{tabular}{lcccccr} % alignment of each column data
		\toprule[\heavyrulewidth]\toprule[\heavyrulewidth]
		\textbf{$n$ }  & \textbf{P} & \textbf{MLE} & \textbf{PMLE} &  \textbf{AIC } \\ 
		\midrule
		\multirow{3}{*}{$202$} & $\alpha$ & $0.02209$ & $0.02208$ & \multirow{3}{*}{$-2810.513$} \\
		& $\beta$ &  $0.00761$ & $0.00761$ \\
		& $\gamma$ & $0.00000$ & $0.00000$\\
			\bottomrule[\heavyrulewidth] 
		\bigskip
			\end{tabular}
\end{table}

\begin{table}[ht!]
	\caption{Model II on Body Mass Index ($X$) and Lean Body Mass ($Y$) }  
	\label{aisI}
	\small % text size of table content
	\centering % center the table
	\begin{tabular}{lcccccr} % alignment of each column data
		\toprule[\heavyrulewidth]\toprule[\heavyrulewidth]
		\textbf{$n$ }  & \textbf{P} & \textbf{MLE} & \textbf{PMLE} &  \textbf{AIC } \\ 
		\midrule
		\multirow{2}{*}{$202$} & $\alpha$ & $0.02210$ & $0.02208$ & \multirow{2}{*}{$-2812.513$} \\
		& $\beta$ &  $0.00761$ & $0.00761$ \\
		
		\bottomrule[\heavyrulewidth] 
		\bigskip
		
	\end{tabular}
\end{table}

\begin{table}[ht!]
	\caption{Model III on  Body Mass Index ($X$) and Lean Body Mass ($Y$) }  
	\label{aisBI}
	\small % text size of table content
	\centering % center the table
	\begin{tabular}{lcccccr} % alignment of each column data
		\toprule[\heavyrulewidth]\toprule[\heavyrulewidth]
		\textbf{$n$ }  & \textbf{P} & \textbf{MLE} &  \textbf{AIC } \\ 
		\midrule
		\multirow{4}{*}{$202$} & $\mu_1$ &  $22.956$ & \multirow{4}{*}{$-2952.242$} \\
		& $\mu_2$ & $64.874$ \\
		& $\sigma^2_1$  & $89.073$\\
		& $\sigma^2_2$  & $170.830$\\
		
		\bottomrule[\heavyrulewidth] 
		\bigskip
		
	\end{tabular}
\end{table}

Note that using the AIC criterion, for the Australian Institute of Sport data on  Mass Index ($X$) and Lean Body Mass ($Y$), we recommend the bivariate dependent equi-dispersed normal conditional model.

\section{Normal variables with variance equal to mean 
		squared, univariate and bivariate}
	
	Consider a normally distributed random variable $X$ with its variance equal to the square of its mean, i.e.,  $X\sim Normal(\tau,\tau^2)$ for some $\tau \in (-\infty,\infty).$ The density of such a random variable $X$ is of the form
	\begin{eqnarray} \label{norm-var-meansq}
	f_X(x;\tau)&=&\frac{1}{|\tau|\sqrt{2\pi}}e^{-(x-\tau)^2/2\tau^2} \nonumber \\
	&& \\
	&=&\frac{1}{|\tau|\sqrt{2\pi}}e^{(x-0.5)} exp{\left\{\frac{-x^2}{                                                                                                                                                                                 2\tau^2}+\frac{x}{\tau}     \right\}}, \nonumber
	\end{eqnarray}
	This  is a curved exponential family and consequently we cannot utilize the Arnold-Strauss theorem to identify the class of all bivariate densities with conditionals in this family. However, since we will be deaing with normal conditionals, we will have conditional moments of the forms displayed in (\ref{e3.100})-(\ref{e3.103}) . From these  equations it is shown in  Appendix A that, in order to have conditional variances equal to the squares of conditional means we will require that
	$$a_{11}=a_{12}=a_{21}=a_{22}=0$$
	and that $a_{20}=a_{10}^2 /2$ and $a_{02}=a_{01}^2/2.$ In such a case, $X$ and $Y$ wil be independent normal variables with variances equal to their means squared. Consequently, the family of bivariate densities  with normal conditionals and with conditional variances equal to squared conditional means is too restrictive to be of interest or of utility.

	Instead, if we think that the class of normal densities with variance equal to the mean squared will be useful to model either marginal or conditional aspects of our data, there are two avenues open to us. First we may consider $(X,Y)$ to have a classical bivariate normal distribution but with two restrictions on the parameters to ensure that $var(X)=[E(X)]^2$ and $var(Y)=[E(Y)]^2$. Such distributions will have marginals in the family (\ref{norm-var-meansq}) but will only have conditionals in that family in the case in which $X$ and $Y$ are independent.
	
	A second approach utilizes the concept of pseudo distributions as described in Filus, Filus and Arnold (2009).
	In this set-uo we postulate that $X$ has a density in the normal with var=mean-squared class, i.e., with desisity of the form  (\ref{norm-var-meansq}) and that for each $x$ the conditional density of $Y$ given $X=x$ is in the class  (\ref{norm-var-meansq}) with a parameter $\tau$ that can depend on $x$. The corresponding joint density will be of the form
	\begin{equation} \label{pseudo-density}
	f_{X,Y}(x,y)=\frac{1}{|\tau_1|\sqrt{2\pi}}e^{(x-0.5)} exp{\left\{\frac{-x^2}{                                                                                                                                                                                 2\tau_1^2}+\frac{x}{\tau_1}     \right\}}\frac{1}{|\tau(x)|\sqrt{2\pi}}e^{(y-0.5)}  exp{\left\{\frac{-y^2}{                                                                                                                                                                                 2\tau(x)^2}+\frac{y}{\tau(x)}     \right\}},
	\end{equation}
	where $\tau_1 \in (-\infty,\infty)$ and $\tau(x)$ is a real valued function. Typically $\tau(x)$ is taken to be a relatively simple function depending on a small number of parameters. For example, we could set $\tau(x)=\tau_2+\tau_3x$ to yield a three-parameter family of denities with the marginal density of $X$ in the class  (\ref{norm-var-meansq}) and all conditional desities of $Y$ given $X=x$ also in the class  (\ref{norm-var-meansq}). A parallel competing model will be one in which the roles of $X$ and $Y$ are interchanged. In practice it will often be difficult to know in advance which of the two models will best fit a given data set and both might be investigated.

	\section{Final remarks} 
	The univariate equi-dispersed normal model was used to construct a corresponding conditionally specified bivariate distribution.  This flexible bivariate model can exhibit a variety of distributional properties including asymmetry, multimodality, marginal skewness and a range of dependence qualities including independence as a special case. A simulation sudy and application to two well-known real data sets, indicate the feasibility of parametric inference for this model. For the two data sets that were considered, the bivariate equi-dispersed normal conditional model provided a better fit than the competing models that were considered. Because the model is flexible even though relatively simply described, it is suggested that it will be a useful addition to the toolkit  of modellers dealing with data that exhibits  skewness, multi-modality and diverse dependence structure.

		\section{Acknowledgment(s)}
	
	The second author’s research was sponsored by the Institution of Eminence (IoE), University of Hyderabad (UoH-IoE-RC2-21-013).

\section{Appendix A}
	
	If we wish to consider bivariate densitites with normal conditionals that will have conditional variances equal to the squares of the corresponding conditional means, we will not be able to find many such distributions. One class of solutions are those which have independent normal marginals with variances equal to the squares of their means. We claim that this is the only valid solution.
	To see this, we may argue as fokllows. 
	
	First observe that since we will have normal conditionals, the conditional means and variances will be,
	as we saw earlier, given by
	\begin{eqnarray*}
	E(X\mid Y=y) & =& -{{a_{12}y^2+a_{11}y+a_{10}} \over
		{2(a_{22}y^2+a_{21}y+a_{20})}} \\ 
	var(X\mid Y=y) & =& {1 \over {2(a_{22}y^2+a_{21}y+a_{20})}}\\
	E(Y\mid X=x) &  =&-{{a_{21}x^2+a_{11}x+a_{01}} \over
		{2(a_{22}x^2+a_{12}x+a_{02})}} \\
	var(Y\mid X=x) & = & {1 \over {2(a_{22}x^2+a_{12}x+a_{02})}}
	\end{eqnarray*}
	If the condition $var(X|Y=y)=[E(X|Y=y)]^2$ is to hold for every $y$, we must then have
	$${1 \over {2(a_{22}y^2+a_{21}y+a_{20})}}=\left[-{{a_{12}y^2+a_{11}y+a_{10}} \over
		{2(a_{22}y^2+a_{21}y+a_{20})}}\right]^2.$$
	Equivalently it must be true that
	$$2(a_{22}y^2+a_{21}y+a_{20})-[a_{12}y^2+a_{11}y+a_{10}]^2=0$$
	for every $y$. The left side is a polynomial of degree 4 and for it to be equal to $0$ for every $y$ , all of its coefficients must be equal to $0$. This implies the following relations must hold
	$$a_{12}^2=0,  \ \ \ \  \  2a_{12}a_{11}=0, \ \ \ \ \ a_{11}^2+2a_{12}a_{10}-2a_{22}=0, \ \ \ \ \ 2a_{11}a_{10}-2a_{21}=0, \ \ \ \ \ a_{10}^2-2a_{20}.$$
	In parallel fashion, if the condition $var(Y|X=x)=[E(Y|X=x)]^2$ is to hold for every $x$, we must then have
	$$a_{21}^2=0,  \ \ \ \  \  2a_{21}a_{11}=0, \ \ \ \ \ a_{11}^2+2a_{21}a_{01}-2a_{22}=0, \ \ \ \ \ 2a_{11}a_{01}-2a_{12}=0, \ \ \ \ \ a_{01}^2-2a_{02}.$$
	
	Clearly we must have $a_{12}=a_{21}=0$, and upon substituting these values our conditions for conditional variances to be equal to squared conditional means simplify to become:
	\begin{equation}\label{keycondition1}
	2(a_{22}y^2+a_{20})=[a_{11}y+a_{10}]^2
	\end{equation}
	and
	\begin{equation}\label{keycondition2}
	2(a_{22}y^2+a_{02})=[a_{11}y+a_{01}]^2.
	\end{equation}
	
	Now, if $a_{11}=0$ then, necessarily, from either equation, $a_{22}=0$. In this case we have constant conditional variances and the model reduces to become a classical bivariate normal one. Moreover, in this case the joint density will factor and thus $X$ and $Y$ are independent, with now the marginal variances equal to the squared marginal means. 
	
	However we must also consider the case in which $a_{11} \neq 0$.  If this is true, then from (\ref{keycondition1}) and (\ref{keycondition2}) it follows that $a_{10}=a_{01}=0$. Then it follows, using the same equations, that $a_{20}=a_{02}=0$.
	But then (\ref{e3.81}) cannot be satisfied and the model is not admissible as a normal conditionals density (it will fail to be integrable). Thus we have confirmed that the only solution has independent normal marginals with variances equal to the squares of their means. 
	
	\section{Appendix B}
	
	Instead of seeking normal-conditionals densities with equi-dispersed conditional densities, we may consider the class of all normal-conditionals densities whose conditional means uniformly exceed the corresponding conditional variances.
	
	If $(X,Y)$ has normal conditionals then its joint density will be of the form (\ref{e3.26}) with 
	conditional moments (once more) of the form
	\begin{eqnarray*}
	E(X\mid Y=y) & = & \mu _{1}(y)=-{{a_{12}y^2+a_{11}y+a_{10}} \over
		{2(a_{22}y^2+a_{21}y+a_{20})}} \\ 
	var(X\mid Y=y) & = & \sigma_{1}^2(y)={1 \over {2(a_{22}y^2+a_{21}y+a_{20})}}\\
	E(Y\mid X=x) & = & \mu _{2}(x)=-{{a_{21}x^2+a_{11}x+a_{01}} \over
		{2(a_{22}x^2+a_{12}x+a_{02})}}\\
	var(Y\mid X=x) & = & \sigma_{2}^2(x)={1 \over {2(a_{22}x^2+a_{12}x+a_{02})}}
	\end{eqnarray*}
	
	If $a_{22}=0$ then also $a_{21}=a_{12}=0$ and the model must be classical bivariate normal. 
	In this case conditional means are linear functions and conditional variances are 
	constants.

	The only examples in this class with conditional means exceeding conditional variances are ones with 
	independent normal marginals.

	If $a_{22}>0$ then there are two constraints on the $a_{ij}$'s in order to have positive conditional variances. They are

	\begin{eqnarray}
	a_{12}^2<4a_{22}a_{02}, \label{a1}\\
	\nonumber \\
	a_{21}^2<4a_{22}a_{20} .\label{a2}
	\end{eqnarray}
	
	In order to have conditional means exceeding conditional variances, the following two quadratic equations must have no real roots:
	\begin{eqnarray}
	a_{12}y^2+a_{11}y+a_{10}+1=0  ,\\
	\nonumber \\
	a_{21}x^2+a_{11}x+a_{01}+1=0.
	\end{eqnarray}
	For this to be true the $a_{ij}$'s must satisfy the following two additional constraints:
	
	\begin{eqnarray}
	4a_{12}(a_{10}+1)>a_{11}^2,  \label{a3}\\
	\nonumber \\
	4a_{21}(a_{01}+1)>a_{11}^2, \label{a4}
	\end{eqnarray}
	
	The $a_{ij}$'s must thus satisfy the four conditions (\ref{a1}),(\ref{a2}),(\ref{a3}) and(\ref{a4}). In addition we must have $a_{12}<0$ and $a_{21}<0$ There are many solutions. For a simple example, 
	set $a_{12}= a_{21}=-1, a_{10}=a_{01}=-2,$ and  $a_{22}=a_{20}=a_{02}=1.$
	
	If, instead, we wish to identify normal-conditionals densities whose conditional means are uniformly less than the corresponding conditional variances, we must impose the same four conditions  (\ref{a1}),(\ref{a2}),(\ref{a3}) and(\ref{a4}), but this time , in addition, we must have $a_{12}>0$ and $a_{21}>0$. There are many solutions in this case also.
	
	\section{Appendix C}
	R code for maximizing nested likelihood function (using closure)  to compute the maximum likelihood estimates for the bivariate equi-dispersed normal conditional model.
	
	\begin{lstlisting}[language=R]
# sample observations or data set
data = data.frame(x,y)	

# nested likelihood function		
my_min = function(data,am,bm,gm)
{
  # Sampling Size			
  nm = dim(data)[1]
		
  # Computing normalizing constant for each choices of am,bm and gm			
  fnm = function(X)
  {
	exp(-(am* X[1]^2 + bm* X[2]^2 + 
	gm * X[1]^2 * X[2]^2 - X[1] - X[2]))
  }
  x = data$x
  y = data$y
  IIn = hcubature(fnm,lower = c(-Inf, -Inf), upper = c(Inf, Inf))$integral    
  lm = -nm * log(IIn) - am * sum(x^2) - bm * sum(y^2) 
  - gm * sum(x^2 * y^2) + sum(x) + sum(y)
  # return -log likelihood
  return(-lm)			
}

# closure function for optimizing 	nested likelihood	
mylik <-function(par)
{
   my_min(data,par[1],par[2],par[3])
}		

# Initial values for am, bm and gm are from optimizing 
# pseudo-likelihood function, say pea, peb and peg	

aI = round(pea,4)
bI = round(peb,4)
gI = round(peg,4)

# maximizing likelihood function: library(optimr)  method: nlminb
result =  optimr(par=c(aI,bI,gI),lower =c(0,0,0),fn=mylik,
control=list(maxit=1, trace=0), method="nlminb")	
\end{lstlisting}

\end{document}